\documentclass[a4paper,11pt, 
]{article}

\usepackage[english]{babel}
\usepackage[latin1]{inputenc}
\usepackage{graphicx}
\usepackage{amscd}
\usepackage{amsfonts}
\usepackage{latexsym}
\usepackage{amsthm}
\usepackage{amssymb,amsmath}
\usepackage{enumerate}
\usepackage{verbatim}

\usepackage{import}

\usepackage[titletoc,title]{appendix}

\usepackage[usenames]{color}

\usepackage{hyperref}
\hypersetup{colorlinks=false}

\usepackage{subfigure}

\title{\bf Tongues in Degree 4 {B}laschke Products}

\author{\small Jordi Canela \thanks{ The three authors were supported by the Spanish project MTM2011-26995-C02-02 and the Catalan project CIRIT 2009-SGR792. The first author was also supported by the Spanish government grant FPU AP2009-4564 and the  grant  346300 for IMPAN from the Simons Foundation.} \\
{\small Inst. of Math. Polish Academy of Sciences (IMPAN)}\\
{\small ul. \'Sniadeckich 8}\\
{\small 00-656 Warszawa, Poland} \and 
{\small N\'uria Fagella $^*$}\\
{\small Dep.~de Matem\`atica Aplicada i An\`alisi}\\
{\small Universitat de Barcelona}\\
{\small Gran Via de les Corts Catalanes, 585}\\
{\small 08005 Barcelona, Spain} \and 
{\small Antonio Garijo $^*$}\\
{\small Dept. d'Enginyeria Inform\`atica i Matem\`atiques}\\
{\small Universitat Rovira i Virgili}\\
{\small Av. Pa\"isos Catalans 26}\\
{\small Tarragona 43007, Spain}  }

\newtheorem{teor}{Theorem} [section]
\newtheorem{thm}[teor]{Theorem}
\newtheorem{propo}[teor]{Proposition}

\newtheorem{lemma}[teor]{Lemma}
\newtheorem{co}[teor]{Corollary}
\newtheorem{conj}[teor]{Conjecture}
\newtheorem*{teoremB}{Theorem B}
\newtheorem*{teoremA}{Theorem A}
\newtheorem*{teoremC}{Theorem C}

\theoremstyle{definition}
\newtheorem{defin}[teor]{Definition}
\newtheorem{defi}[teor]{Definition}
\newtheorem{rem}[teor]{Remark}

\newsavebox{\savepar}

\newcommand{\com}{\mathbb{C}}
\newcommand{\wcom}{\widehat{\mathbb{C}}}

\newcommand{\real}{\mathbb{R}}

\newcommand{\nat}{\mathbb{N}}

\newcommand{\dis}{\mathbb{D}}
\newcommand{\cercle}{\mathbb{S}^1}
\newcommand{\re}{\rm{Re}}
\newcommand{\im}{\rm{Im}}

\newcommand{\modul}{\rm{mod}\;}
\newcommand{\XX}{\mathcal{X}}
\newcommand{\inv}{\mathcal{I}}

\date{\today}

    {\end{list}}

\begin{document}
\maketitle

\begin{abstract}
 {\noindent \small The goal of this paper is to investigate  the family of Blasche products  $B_a(z)=z^3\frac{z-a}{1-\bar{a}z}$, which is a rational family of perturbations of the doubling map. We focus on the tongue-like sets which appear in its parameter plane. We first study their basic topological properties and afterwords we investigate how bifurcations take place in a neighborhood of their tips. Finally we see how the  period one tongue extends beyond its natural domain of definition.
 }

\vspace{0.5cm}
\textit{Keywords: holomorphic dynamics, Blaschke products, circle maps, tongues.}
\end{abstract}


\section{Introduction}

 Given a rational map $f:\widehat{\com}\rightarrow \widehat{\com}$, where $\widehat{\com}=\com\cup \{\infty\}$ denotes  the  Riemann sphere, we consider the dynamical system given by the iterates of $f$. The Riemann sphere splits into two totally $f$-invariant subsets: the  \textit{Fatou set} $\mathcal{F}(f)$, which  is defined to be the set of points $z\in\widehat{\com}$ where the family  $\{f^n, n\in\mathbb{N}\}$ is normal in some neighborhood of $z$, and  its complement,  the \textit{Julia set} $\mathcal{J}(f)$. The dynamics of the points in $\mathcal{F}(f)$ are stable in the sense of normality or equicontinuity whereas the dynamics in $\mathcal{J}(f)$ present chaotic behavior. The Fatou set is open and its connected components, called \textit{Fatou components}, are mapped under $f$ among themselves. All Fatou components of a rational map are either periodic or preperiodic (see \cite{Su}). By means of the Classification Theorem (c.f.~\cite{Mi1}), any periodic Fatou component of a rational map  is either the basin of attraction of an attracting or parabolic cycle, or  a simply connected rotation domain (a Siegel disk), or a doubly connected rotation domain (a Herman ring). Moreover, any such component is somehow related to a \textit{critical point}, i.e.\ a point $z\in\widehat{\com}$ such that $f'(z)=0$. Indeed, the basin of attraction of an attracting or parabolic cycle contains, at least, a critical point whereas Siegel disks and Herman rings have critical orbits (orbits of critical points) accumulating on their boundaries.

The aim of this paper is to study the bifurcations that occur in the parameter plane of the degree 4 Blaschke products given by 

\begin{equation}\label{blasformula}
B_{a}(z)=z^3\frac{z-a}{1-\bar{a}z},
\end{equation}

\noindent where $a,z\in\com$. As we shall explain below, this family is interesting for several reasons, being the rational analogue of the well known double standard family of perturbations of the doubling map of the unit circle. Indeed, as all finite Blaschke products, the maps $B_a$ leave the unit circle $\cercle$ invariant, i.e.\  $B_a(\cercle)=\cercle$.  Moreover, they are rational perturbations of the doubling map of the circle $R_2(z)=z^2$ (equivalently given by $\theta\rightarrow 2\theta \; (\modul\; 1))$. Indeed, if $a=re^{2\pi i\alpha}$ with $\alpha\in\real$, the rational maps $B_a$ converge uniformly over compact subsets of the punctured plane $\com^*$ to $e^{4\pi i \alpha}z^2$ as $r$ tends to $\infty$. In fact, if $|a|\geq2$, the circle maps $B_a|_{\cercle}$ are degree two coverings of the unit circle and hence they are semiconjugate to the doubling map (c.f.~\cite{MiRo1}).

For all values of $a\in\com$, the points $z=0$ and $z=\infty$ are superattracting fixed points of local degree 3 (c.f.\ \cite{CFG1}). We denote by $A_a(0)$ and $A_a(\infty)$ their basins of attraction and by $A_a^*(0)$ and $A_a^*(\infty)$ their immediate basins of attraction, i.e.\ the connected components of $A_a(0)$ and $A_a(\infty)$ which contain $z=0$ and $z=\infty$, respectively. We shall drop the dependence on $a$ whenever it is clear from the context. If $|a|\leq 1$, $B_a(\dis)=\dis$ and the basins of attraction $A(0)=\dis$ and $A(\infty)=\wcom\setminus\overline{\dis}$ are the only Fatou components, separated by the Julia set which is necessarily $\cercle$. But for every other parameter, there is a preimage of infinity $z_{\infty}\in\dis$ and the two free critical points $c_{\pm}$ (distinct unless $|a|=2$)  may lead to the existence of stable components different from $A(0)$ and $A(\infty)$. If $|a|\geq2$, even if there are two free critical points, the Blaschke family is essentially unicritical due to the symmetry with respect to $\cercle$ which, in this case, ties the two critical orbits together in a certain sense.

Given a family $f_a$ of orientation preserving homeomorphisms of the unit circle, where $a\in\Delta$ and $\Delta\subset\wcom $ or $\Delta\subset \real^2$, one can assign a rotation number to each of its members which gives the average asymptotic rate of rotation of points in the circle. The level sets of the rotation number are known as tongues or Arnold tongues, since they were introduced by  V.~Arnold \cite{Ar} for the standard family of perturbations of the rigid rotation

$$\theta\rightarrow \theta+\alpha+(\beta/2\pi)\sin(2\pi\theta) \; \; (\mbox{mod}\; 1)$$

\noindent where $0\leq\theta<1$, $0\leq \alpha\leq 1$ and $0\leq \beta\leq1$. Tongues have been studied by many authors, both for the Arnold standard maps (see e.g.\ \cite{Boyl}, \cite{WBJ}, \cite{EKT} or \cite{dLL}) and for other families of homeomorphisms of the unit circle such as the degree~$3$ Blaschke products introduced by Herman \cite{Her}. Of special relevance are the Arnold tongues corresponding to rational rotation numbers, since those are precisely the parameters for which the maps have attracting or parabolic periodic orbits on $\cercle$. If the maps $f_a$ are not homeomorphisms but degree $2$ covers of $\cercle$, we cannot assign a rotation number to them. However, tongues can still be defined as sets of parameters for which $f_a$ has an attracting cycle in $\cercle$. We can associate a type $\tau(a)$ to every such $f_a$, where $\tau(a)$ is a periodic point of the doubling map and hence describes the rotation patterns of the attracting cycle of $f_a$ (see Section~\ref{stongues} for details). In this setting, a tongue $T_{\tau}$ is defined as the open set of parameters $a\in\Delta$ of type $\tau$. They were studied by  M.~Misiurewicz and A.~Rodrigues \cite{MiRo1,MiRo2} for the double standard family of perturbations of the doubling map

$$\theta\rightarrow 2\theta+\alpha+(\beta/\pi)\sin(2\pi\theta) \;\; (\mbox{mod} \; 1)$$

\noindent where $0\leq\theta<1$, $0\leq \alpha\leq 1$ and $0\leq \beta\leq1$.  Later on, A. Dezotti \cite{De} used the complex extension of the double standard maps on the punctured plane, which is given by  $z\rightarrow e^{i\alpha}z^2e^{\beta/2(z-1/z)}$, in order to prove the connectivity of the tongues. This family was also studied by R.~de~la~Llave, M.~Shub and C. Sim\'o \cite{dLSS}, who considered entropy related issues of the $k$-th standard maps $\theta\rightarrow k \theta+\alpha+\epsilon\sin(2\pi\theta)  \; (\mbox{mod}\; 1)$ for $\epsilon$ small and $k\geq2$.

Given that the Blaschke products $B_a$ are rational perturbations of the doubling map, $B_a|_{\cercle}$ may be considered as the rational analogue of the double standard family. Although there is no explicit simple expression for the restriction of $B_a$ to $\cercle$, the global dynamics are simpler than in the transcendental case. If $|a|\geq2$, the $B_a|_{\cercle}$ are degree $2$ coverings of the unit circle and the tongues are well defined. Inspired by the mentioned works of Misiurewicz, Rodriguez and Dezotti, we show that they  are connected and simply connected. More precisely, we prove the following theorem.

\begin{teoremA}
Given any periodic point $\tau$ of the doubling map the following results hold.

\begin{enumerate}[(a)]
\item  The tongue $T_{\tau}$ is not empty and consists of three connected  components (only one connected component if we consider the parameter plane modulo the symmetries given by the third roots of the unity).
\item Each connected component of $T_{\tau}$ contains a unique parameter $r_{\tau}$, called the \emph{root of the tongue}, such that $B_{r_{\tau}}$ has a superattracting cycle in $\cercle$. The root $r_{\tau}$ satisfies $|r_{\tau}|=2$. 
\item Every connected component of $T_{\tau}$ is simply connected.

\item The boundary of every connected component of $ T_{\tau}$  consists of two curves which are continuous graphs as function of $|a|$ and intersect each other in a unique parameter $a_{\tau}$ called the \emph{tip of the tongue}.

\end{enumerate}

\end{teoremA}

The boundary of a tongue $T_{\tau}$ of period $p$ is the union of two curves of parameters which intersect at the tip $a_{\tau}$ of the tongue (see Figure~\ref{zoomtip}). These boundary parameters correspond to maps $B_a$ which have a persistent parabolic cycle of period $p$ and multiplier $1$ in $\cercle$. On the tip $a_{\tau}$ of the tongue, the parabolic cycle of $B_{a_{\tau}}$ has multiplicity 3, i.e.\ $B_{a_{\tau}}$ has a cycle $\{z_0,\cdots,z_{p-1}\}$ of period $p$,  multiplier $1$ and second derivative zero ($(B_{a_{\tau}}^p)''(z_0)=0$) in $\cercle$. Along $\partial T_{\tau}\setminus a_{\tau}$ there is a persistent saddle-node bifurcation taking place: two period $p$ cycles collide in $\cercle$ and exit it (see Figure~\ref{bifneartip}).  Notice that these bifurcation along curves cannot happen for a uniparametric family of holomorphic maps which depend holomorphically on the parameter. The real saddle-node bifurcation was studied by M.~Misiurewicz and R.~A.~P\'erez \cite{MiP}  from a complex point of view. They characterized, depending on the sign of the Schwarzian derivative, whether the period $p$ cycles exiting the unit circle (or the real line) are attracting or repelling.  Crowe et al \cite{Cr} showed that this type of bifurcations also occurs in the Tricorn, the connectedness locus of the antipolynomials $p_c(z)=\overline{z}^2+c$. Their result was later generalized by  J.~H.~Hubbard and D.~Schleicher \cite{HS}. They studied these bifurcations in the Multicorns, the bifurcation loci of the antipolynomials $p_{d,a}(z)=\overline{z}^d+a$, by using the holomorphic index of the fixed points. Using the holomorphic index and results of algebraic geometry such as Chevalley's Theorem (see Theorem~\ref{Chevalley}, c.f.~\cite{Ha}), we prove the following result.

\begin{teoremB}
Let $a_{\tau}$ be the tip of a tongue $T_{\tau}$ of period $p$. Then, there exists a neighborhood $U$ of $a_{\tau}$ such that if $a\in U$ then, either $a\in T_{\tau}$ or $a\in \partial T_{\tau}$ or $a$ belongs to a  hyperbolic component of disjoint type, i.e.\ a hyperbolic component for which the maps $B_a$ have two attracting cycles other than $z=0$ and $z=\infty$.
\end{teoremB}

Using the same techniques utilized in Theorem~B and the parametrization of disjoint hyperbolic components described in \cite{CFG1} (see Theorem~\ref{thmparametrize}) we also proof that any parameter $a$ with $|a|>2$ such that $B_a$ has a parabolic cycle of multiplicity 3 in $\cercle$ corresponds to a tip of a tongue (see Proposition~\ref{noparsol}). We conclude that any parameter $a$ with $|a|>2$ such that $B_a$ has a parabolic cycle in $\cercle$ belongs to the boundary of a tongue (see Corollary~\ref{parabolicinboundary}).

If $a$ belongs to the open annulus $\mathbb{A}_{1,2}$  of inner radius $1$ and outer radius $2$, then the Blaschke products $B_a|_{\cercle}$ are no longer degree two coverings of the unit circle. Despite that,  the attracting cycle associated to a given tongue may be continued for parameters within $\mathbb{A}_{1,2}$. This leads to the concept of \emph{extended tongues} $ET_{\tau}$, open connected sets of parameters $a$, $|a|>1$, for which $B_a|_{\cercle}$ has an attracting cycle which can be real analytically continued to the attracting cycle of a tongue $T_{\tau}$. However, extended tongues are not disjoint. Indeed,  if $a\in \mathbb{A}_{1,2}$, then the two free critical points of $B_a$ lie on the unit circle and their orbits are not related by symmetry. Therefore, they may accumulate on different attracting cycles, allowing a parameter to belong to two different extended tongues simultaneously. We focus on the study of the extended fixed tongue $ET_0$ and prove the following theorem (see Figure~\ref{paramblash3col1}).

\begin{teoremC}
Given two connected components of the fixed tongue $T_0$, the intersection of their extensions in $\mathbb{A}_{1,2}$ is empty. The boundary of every connected component of the extended fixed tongue $ET_0$ consists of two disjoint connected components. The exterior component consists of parameters for which there is a parabolic fixed point of multiplier $1$. The interior component consists of parameters for which there is a parabolic fixed point of multiplier $-1$. Moreover, there is a period doubling bifurcation taking place throughout the curve of interior boundary parameters.
\end{teoremC}

The paper is structured as follows. In Section~\ref{introblas0} we introduce the basic properties of the Blaschke products $B_a$ and present some alternative parametrizations which are useful later on. In Section~\ref{stongues} we introduce the concept of tongues for the family $B_a$ and prove Theorem~A. In Section~\ref{tiptongue} we study the bifurcations which take place around the tip of the tongues proving Theorem~B. Finally, in Section~\ref{exttongue} we see how tongues extend within the annulus $\mathbb{A}_{1,2}$ and study the extended fixed tongue $ET_0$ proving Theorem~C.  

\textit{Acknowledgements.} The authors would like to thank A.~Dezotti, A.~Epstein, M.~Misiurewicz, J.~C.~Naranjo, M.~S\'aez and A.~Vieiro for their many and useful comments which greatly improved this paper. 

\section{Preliminaries on the Blaschke family}\label{introblas0}

In this section we give an introduction to the Blaschke products $B_a$ as in (\ref{blasformula}). We first study the basic properties of the dynamical plane and the parameter plane. Afterwords we consider some reparametrizations of the family which are useful later on.  We refer to \cite{CFG1} for proofs of the results and a more detailed description of the family.

 We consider the degree 4 Blaschke products of the form 
 
\begin{equation}\label{blasformula2}
 B_{a,t}(z)=e^{2\pi it}z^3\frac{z-a}{1-\bar{a}z},
 \end{equation}

\noindent where $a \in \com$ and $t \in \mathbb{R}/\mathbb{Z}$. Since $B_{a,t}$ leaves invariant the unit circle, it is symmetric with respect $\cercle$, i.e.\ $B_{a,t}(z)=\mathcal{I} \circ B_{a,t}\circ \mathcal{I}(z)$ where $\mathcal{I}(z)=1/\bar{z}$.
Even if it is sometimes necessary to consider the parameter $t$, the next lemma tells us that, for the purpose of classification,  we can omit it. The proof is straightforward.

\begin{lemma}\label{conjblas}
Let $\alpha\in \real$ and let $\eta(z)=e^{-2\pi i\alpha}z$. Then $\eta$ conjugates the maps $B_{a,t}$ and $B_{a e^{-2\pi i\alpha}, t+3\alpha}$. In particular, $B_{a,t}$ is conjugate to $B_{a e^{\frac{2\pi i t}{3}}, 0}$.

\end{lemma}

The orbits of the critical points, i.e.\ the points $c\in\wcom$ such that $f'(c)=0$, control the possible stable dynamics of any rational map $f$. Since  the Blaschke products $B_a$  have degree $d=4$, they have $2d-2=6$ critical points counted with multiplicity. The fixed points $z=0$ and $z=\infty$ are critical points of multiplicity 2 and hence superattracting fixed points of local degree 3. The other two  critical points, denoted by $c_{\pm}$, are given by 

\begin{equation}\label{criticalpoints}
c_{\pm}:=c_{\pm}(a):=a \cdot \frac{1}{3|a|^2}\left(2+|a|^2\pm\sqrt{(|a|^2-4)(|a|^2-1)}\right).
\end{equation}

The next lemma states us that the parameters $a$ and $t$ depend continuously on the critical points.

\begin{lemma}[{\cite[Lemma~3.2]{CFG1}}]\label{continuitya}
Given a Blaschke product $B_{a,t}$ as in (\ref{blasformula2}) with $|a|\geq2$ or $|a|< 1$, the parameter $a$ is continuously determined by the critical points $c_{\pm}$. Moreover, if the image $B_{a,t}(z_0)\neq\{0,\infty\}$ of a point $z_0\in\com^*$ is fixed, then $t$ depends continuously on $a$.

\end{lemma}

\begin{figure}[hbt!]
    \centering
    \subfigure[\scriptsize{Case $|a|<1$} ]{
   			 \def\svgwidth{160pt}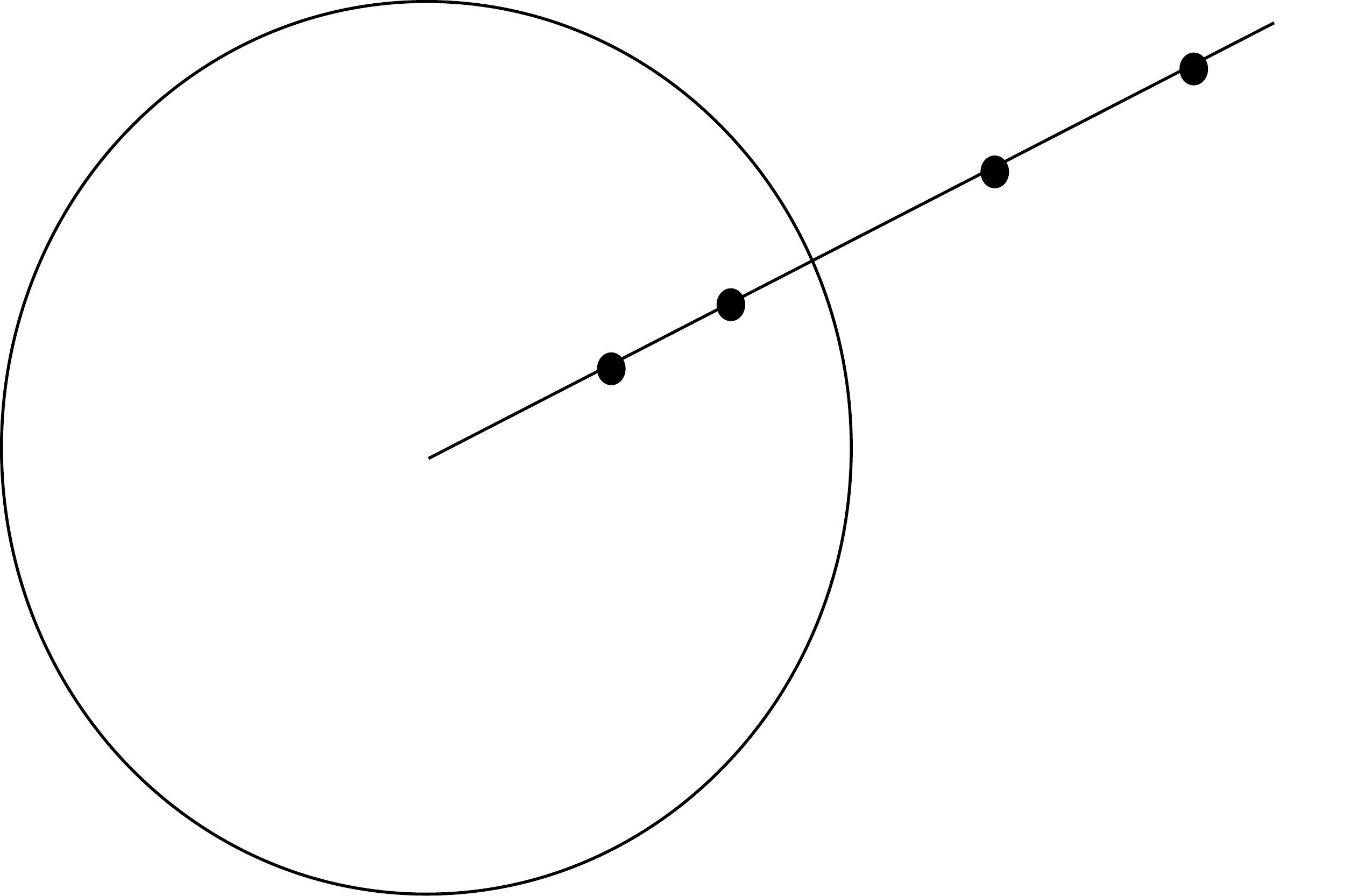
     }
    \subfigure[\scriptsize{Case $1<|a|<2$} ]{
    		 \def\svgwidth{160pt}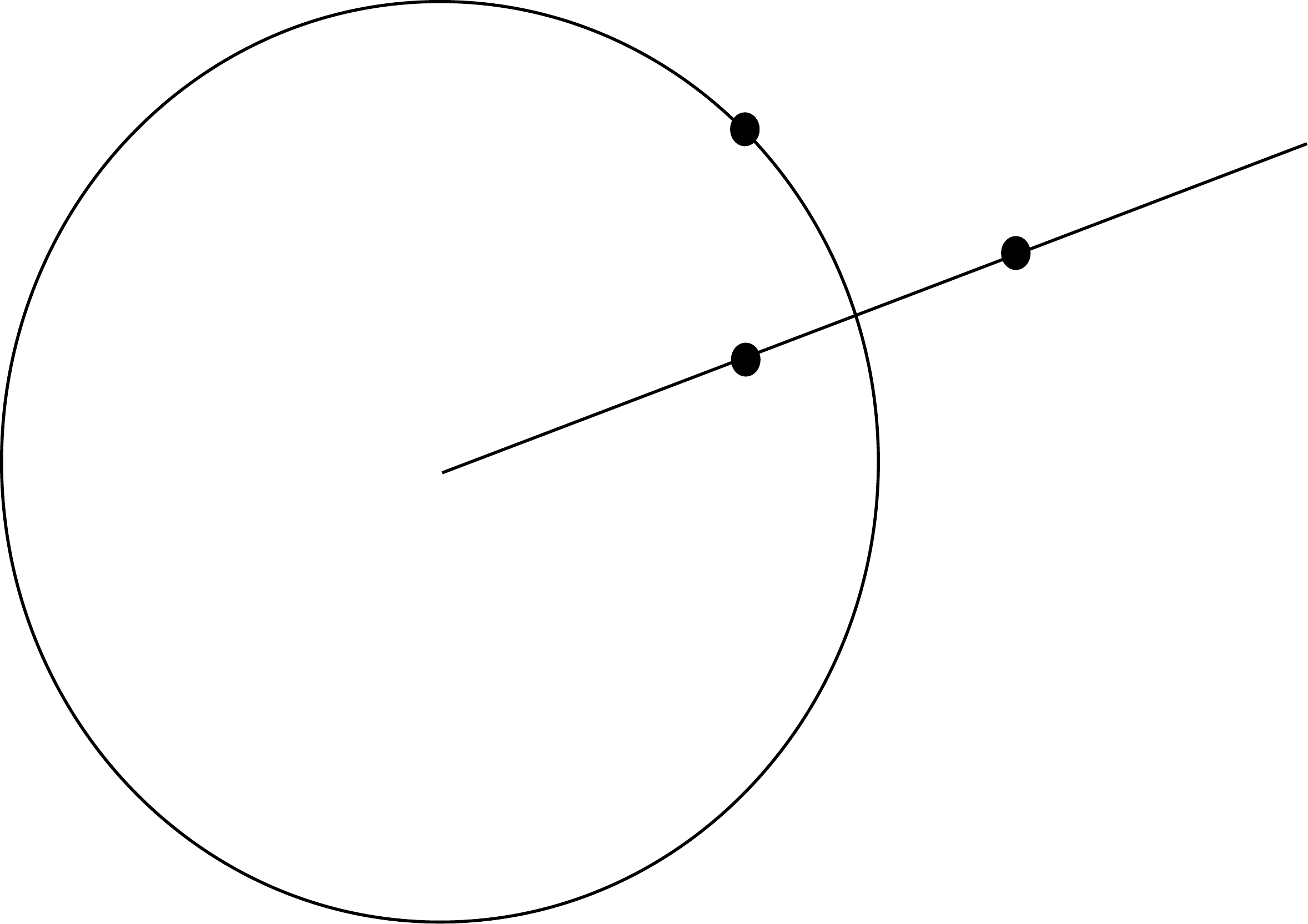}
    \hspace{0.1in}
    \subfigure[\scriptsize{Case $|a|=2$} ]{
 		  \def\svgwidth{160pt}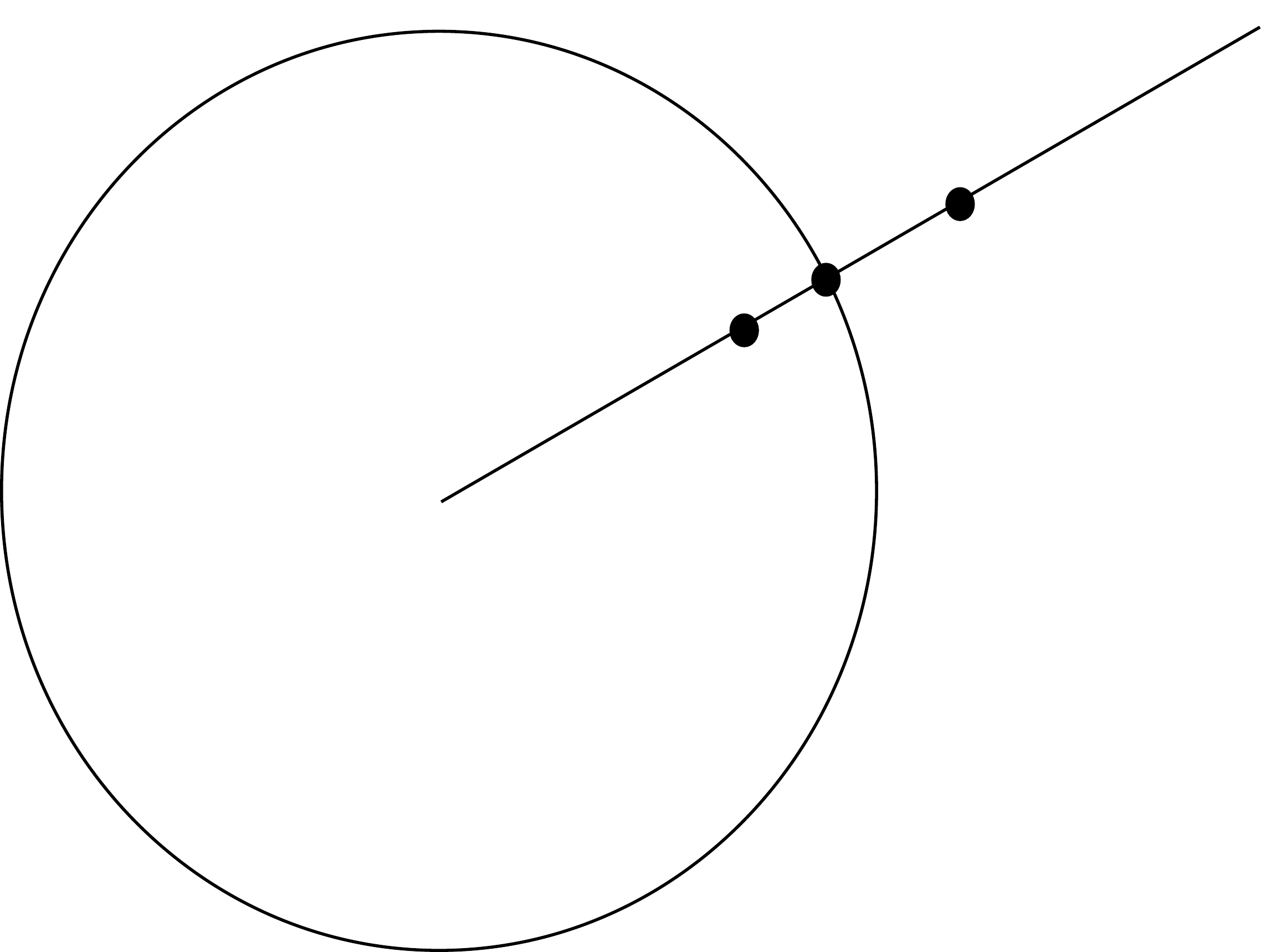}
    \subfigure[\scriptsize{Case $|a|>2$} ]{
   		   \def\svgwidth{160pt}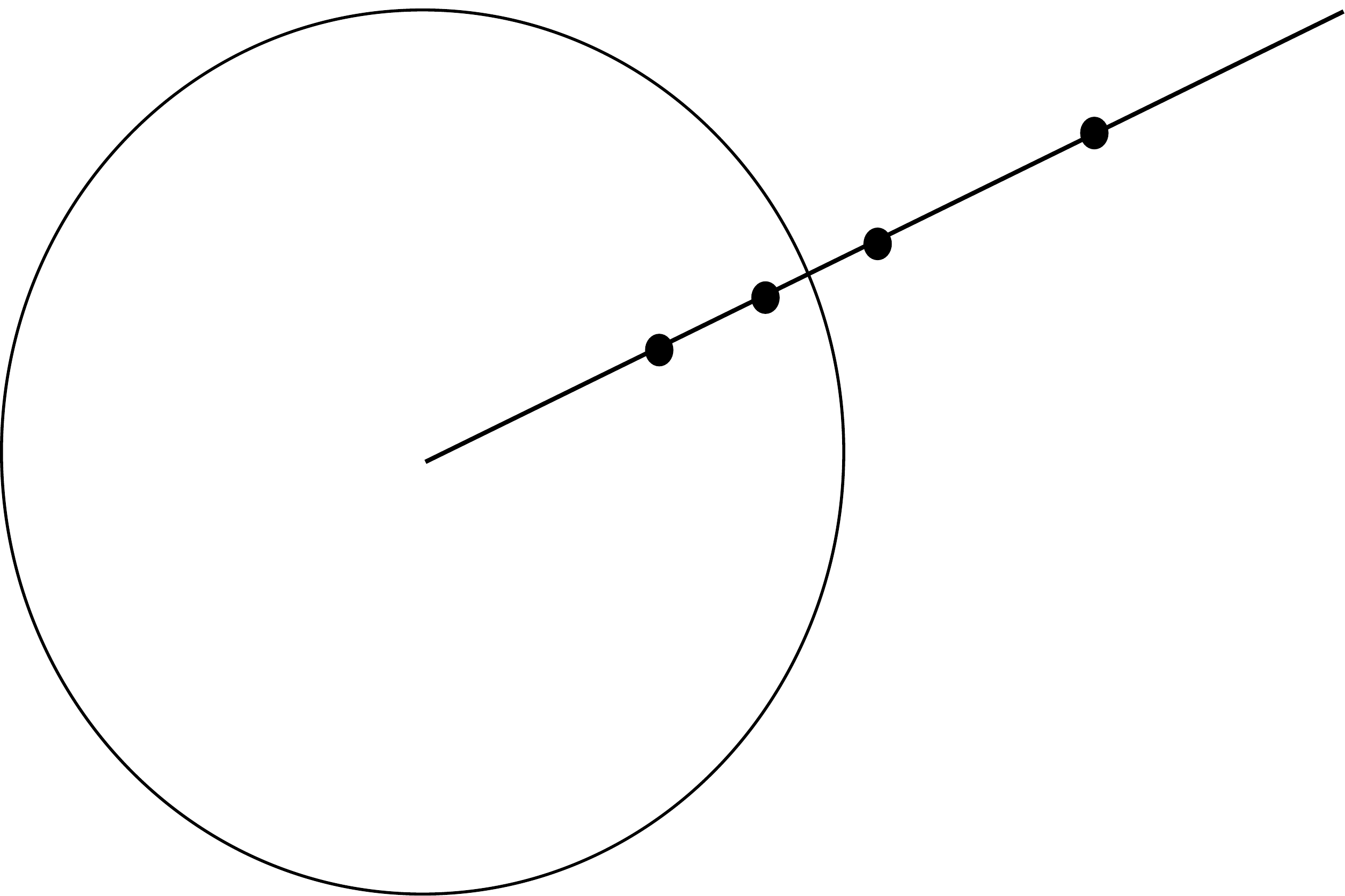}
    \caption{\small{Different configurations of the critical points and the preimages of zero and infinity with respect to the unit circle depending on $|a|$. We also draw the sets of preimages of $\cercle$.}} 
    \label{esquemapunts}
\end{figure}

The positions of the free critical points, together with the ones of the pole $z_\infty=1/\overline{a}$ and the zero $z_0=a$ determine the possible stable dynamics of the Blaschke products. Indeed, if the pole is not in $\dis$, then $B_a(\dis)=\dis$ and all points in $\dis$ tend under iteration of $B_a$ to the superattracting cycle $z=0$. The possible configurations can be summarized, depending on the modulus of $a$, as follows (see Figure~\ref{esquemapunts}).

\begin{itemize}
\item If $|a|<1$ then $z_{\infty}\in\wcom\setminus\overline{\dis}$ and $\mathcal{J}(B_a)=\cercle$.
\item If $|a|=1$ then the two free critical points, the pole and the zero collide in a single point. For these parameters, the Blaschke products degenerate to the family of cubic polynomials $-az^3$ and $\mathcal{J}(B_a)=\cercle$.
\item If $1<|a|<2$ then $z_{\infty}\in\dis$ and $c_{\pm}\in\cercle$. The two critical orbits are independent, which may lead to up to two attracting or parabolic cycles in $\cercle$.
\item If $|a|=2$ then  $z_{\infty}\in\dis$ and the two critical points collide in a single one $c\in\cercle$. There can be at most one attracting cycle in $\cercle$ and $B_a|_{\cercle}$ is a degree two cover.

\item If $|a|>2$ then  $z_{\infty}\in\dis$, the critical points do not belong to $\cercle$ and are symmetric with respect to it, i.e.\ $c_+=1/\overline{c_-}$. Consequently, the critical orbits also are symmetric with respect to $\cercle$ and therefore have the same (or symmetric) asymptotic behaviour. Moreover, $B_a|_{\cercle}$ is a degree two cover.

\end{itemize}

Along this paper we mainly focus on parameters $a$ such that $|a|\geq2$ although we analyse the case $1<|a|<2$ in Section~\ref{exttongue}.

\begin{figure}[hbt!]
\centering
\includegraphics[width= 425pt]{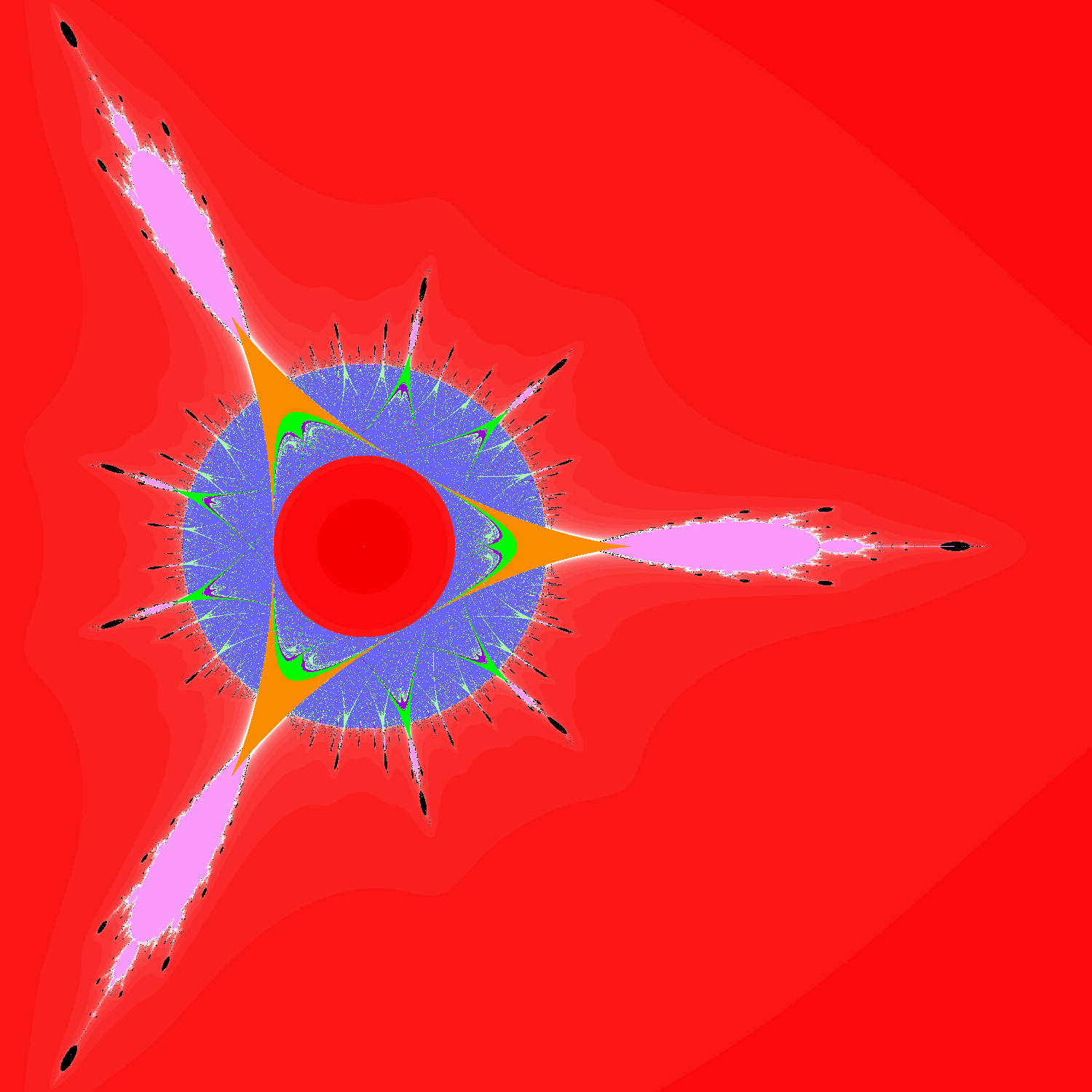}
\caption{\small Parameter plane of the Blaschke family $B_a$. The colors are as follows: red if $c_+\in A(\infty)$, black if  $c_+\in A(0)$; orange, strong green and violet if $O^+(c_+)$ accumulates on a periodic point in $\cercle$ of period 1, 2 and 3, respectively; pallid green if $O^+(c_+)$ accumulates on a periodic point in $\cercle$ with period other than 1, 2 or 4; pink if $O^+(c_+)$ accumulates on a periodic orbit not in $\cercle$ and blue in any other case. The inner red disk corresponds to the unit disk. }
\label{paramblash0}
\end{figure}

We finish this first part of the section stating two results on the parameter plane which will be useful later on. As we can see on Figure~\ref{paramblash0}, there are some symmetries which can be observed in the parameter plane. They are explained in the following lemma (see \cite[Lemma~5.1.1]{Can} for details, c.f.~\cite[Lemma~4.1]{CFG1}).

\begin{lemma}\label{symm}
Let $a, b\in\com\setminus\cercle$. Then $B_a$ and $B_b$ are conformally conjugate if and only if $b=\xi a$ or $b=\xi \overline{a}$, where $\xi$ is a third root of the unity. 
\end{lemma}

The following theorem,  \cite[Theorem~C]{CFG1}, tells us that every disjoint hyperbolic component $\Omega$, i.e.\ every connected set of parameters $a$ such that $B_a$ has two different attracting cycles not contained in $\cercle$ and other than $z=0$ and $z=\infty$, is homeomorphic to the unit disk. The homeomorphism is given by the multiplier map, which associates, to every $a\in\Omega$, the multiplier of the attracting cycle on which the orbit of $c_+$ accumulates. 

\begin{thm}[{\cite[Theorem~C]{CFG1}}]\label{thmparametrize}
Let $\Omega\subset \{a\in\com; |a|>2\}$ be a disjoint hyperbolic component. Then, the multiplier map is a homeomorphism between $\Omega$ and the unit disk.

\end{thm}

\subsection{Reparametrizing the Blaschke family}\label{prelimparametrization}

In this section we describe some alternative parametrizations of $B_{a}$ which are particularly useful when we restrict to the unit circle.  Let $a=re^{2\pi i \alpha}$ with  $r\in (1,\infty)$ and $\alpha\in [0,1/3)$ (or $\alpha\in\real/\frac{1}{3}\mathbb{Z}$). From Lemma~\ref{conjblas} we know  that $B_a=B_{a,0}$ as in (\ref{blasformula}) is conformally conjugate to $B_{r,3\alpha}$ as in (\ref{blasformula2}). It is  enough to restrict to parameters $\alpha\in[0,1/3)$ due to symmetry on the parameter plane (see Lemma~\ref{symm}). Notice that if $\alpha\in[0,1/3)$ we have a one to one correspondence between the parameters $a$ of the family $B_a$ and the parameters $(r,\alpha)$ of $g_{r,\alpha}:=B_{r,3\alpha}|_{\cercle}$. Summarizing, we consider the circle maps
 
\begin{equation}\label{gra}
g_{r,\alpha}(e^{2\pi i x})=e^{6\pi i x}e^{6\pi i \alpha}\frac{e^{2\pi i x}-r}{1-re^{2\pi i x}},
\end{equation}
 
\noindent where $r\in(1,\infty)$ and $\alpha\in [0,1/3)$.  Its  lift has the form

\begin{equation}\label{hra}
h_{r,\alpha}(x)=3x+3\alpha+\frac{1}{2\pi i}\log\left(\frac{e^{2\pi i x}-r}{1-re^{2\pi i x}}\right).
\end{equation}

We shall often use $g_{r,\alpha}$ instead of $B_a|_{\cercle}$ given that its lift is somehow simpler. Indeed, it follows directly from its expression that $h_{r,\alpha}$ is strictly increasing with respect to $\alpha$.

 \begin{lemma}\label{increasing}
Let $r\geq 2$. Then, the lift $h_{r,\alpha}(x)$ satisfies that $\frac{\partial}{\partial x}h_{r,\alpha}(x)$ is non-negative for all $x$. Moreover, for any $p\in \mathbb{N}$, the mapping $\alpha\rightarrow h_{r,\alpha}^{ p}(x)\in \cercle$ is strictly increasing and, if $r\geq3$, then $h_{r,\alpha}^{ p}(x)\geq1$ for all $x, \alpha \in\real$. 

\end{lemma}
\proof

We prove that $\frac{\partial}{\partial x}h_{r,\alpha}(x)$ is non-negative for all $x$, and hence so is $\frac{\partial}{\partial x}h^{ p}_{r,\alpha}(x)$ for all $p$. Then, strict monotonicity with respect to $\alpha$ for all $p$ follows from the fact that we have it for $p=1$. We also prove that $\frac{\partial}{\partial x}h_{r,\alpha}(x)\geq 1$ if $r\geq 3$. Notice that $\frac{\partial}{\partial x}h_{r,\alpha}(x)$ is given by the formula

\begin{equation}\label{derhra}
\frac{\partial}{\partial x} h_{r,\alpha}(x)=3+\frac{1-r^2}{1+r^2-2r \cos(2\pi x) }.
\end{equation}

It can easily be seen that this expression is non-negative for $r\geq 2$. Indeed, the minimum of this function is taken whenever $x=0$, and
$$\frac{\partial}{\partial x}h_{r,\alpha}(0)=3+\frac{1-r^2}{1+r^2-2r }=3+\frac{(1+r)}{(1-r)}.$$
For $r>1$ this is an increasing function which is equal to zero for $r=2$. Moreover, it is  greater than $1$ when $r\geq3$.
\endproof

It will also be useful to consider the circle maps $g_{r,\alpha}$ as restrictions of the rational maps 

\begin{equation}\label{Gab}
G_{a,b}(z)=bz^3\frac{z-a}{1-az},
\end{equation}

\noindent where $a,b\in\com$. The $G_{a,b}$ define a degree $4$ family of almost bicritical maps unless $a=\pm1$, when they degenerate to the degree $3$ polynomials $\mp bz^3$. By almost bicritical family we mean that, as is the case for the Blaschke products $B_a$, the maps $G_{a,b}$ have at most two free critical points whilst the other ones are permanent superattracting fixed points. Notice that these maps are not symmetric with respect to the unit circle. Similarly to the Blaschke products $B_a$, the points $z=0$ and $z=\infty$ are superattracting fixed points of local degree $3$ for all $G_{a,b}$. The two free critical points are given by
$$
c_{\pm}:=c_{\pm}(a):=\frac{1}{3a}\left(2+a^2\pm\sqrt{(a^2-4)(a^2-1)}\right)
$$

\noindent and are the solutions of $3az^2-2(a^2-2)z+3a=0$. In particular we have that $c_+\cdot c_-=1$ and none of them is equal to zero if $a\neq 0$ and $a\neq\infty$.

We say that a parameter $(a,b)$ is  escaping if the orbit of any of the critical points accumulates on $z=0$ or $z=\infty$. We finish this section studying the non-escaping set of the family $G_{a,b}$, i.e.\ the set of parameters $(a,b)$ for which none of the critical orbits accumulates on $z=0$ or $z=\infty$.

\begin{lemma}\label{boundednonescapingGab1}
The non-escaping set of the family $G_{a,b}$ is bounded in the $a$-parameter, i.e.\ there exists a constant $C>0$ such that if $|a|>C$ then the parameter $(a,b)$ is escaping.
\end{lemma}
\proof
 We will prove that, if $|a|$ is big enough then one of the critical orbits accumulates on $z=\infty$ or on $z=0$. Notice that the other critical orbit may accumulate on a bounded attracting cycle even if  $|a|$ tends to infinity.

 We distinguish two cases. Assume first that $|b|\geq 1$.
We have that, if $|z|>\lambda(|a|+1)$ with $\lambda\geq 1$, then $|G_{a,b}(z)|>\lambda|b||z|$. Indeed, one can check that $|z-a|>\lambda$ and that 
$$|z|^2>|z|(|a|+1)=|za|+|z|>|za|+1>|1-az|.$$
Therefore, we have
$$|G_{a,b}(z)|=|b||z|^3\frac{|z-a|}{|1-az|}>|z|^3\frac{\lambda}{|z|^2}=\lambda|z|.$$

 As $|a|$ tends to infinity the critical point $c_+(a)$ tends to $2a/3$ and  $c_-(a)$ tends to $3/(2a)$. Consequently,  the modulus of the critical value $v_+=G_{a,b}(c_+(a))$ grows as $M |a|^2$ for some $M>0$ and, for $|a|$ large enough, $|v_+|>\lambda(|a|+1)$ with $\lambda>1$. We conclude that \linebreak $|G_{a,b}^n(v_+)|\rightarrow\infty$  when $n\rightarrow\infty$. Therefore, if $|a|$ is large enough  and $|b|\geq1$ then the parameter $(a,b)$ is escaping.
 
 Consider now the case $|b|<1$. First we prove that, if $|a|>1$ and $|z|<1/(2|a|)$ then $|G_{a,b}(z)|<3|b||z|/4$. From these  hypothesis we conclude that $|z|<1/2$ and obtain the inequalities
 $$|z-a|<|z|+|a|<\frac{1}{2|a|}+|a|<\frac{|a|}{2}+|a|=3|a|/2$$ 
 and
  $$ |1-az|>1-|az|>1-\frac{1}{2}=\frac{1}{2}.$$
  Therefore, we have
  $$|G_{a,b}(z)|=|b||z|^3\frac{|z-a|}{|1-az|}<|b||z|^3 3|a|<\frac{3}{2}|b||z|^2<\frac{3|b||z|}{4}.$$
  
Since $c_-(a)$ converges to $3/(2a)$ as $a$ tends to infinity, we conclude that the modulus of the critical value $v_-=G_{a,b}(c_-)$ decreases as $M/|a|^2$ for some $M>0$. Hence, for $|a|$ large enough, $|v_-|<1/(2|a|)<1/3$ and $|G_{a,b}^n(v_-)|\rightarrow 0$ when $n\rightarrow\infty$. Therefore, if $|a|$ is large enough  and $|b|<1$ then the parameter $(a,b)$ is escaping.

\endproof

\begin{lemma}\label{boundednonescapingGab2}
For fixed $a_0\in\com$, $a_0\neq \pm 1$, the non-escaping set of $G_{a_0,b}$ is bounded with respect to  the parameter $b$, i.e.\ there exists a constant $C(a_0)>0$ such that if $|b|>C(a_0)$ then the parameter $(a_0,b)$ is escaping.
\end{lemma}
\proof
It is enough to prove that if $|b|$ is large enough then the orbit of one of the critical points accumulates on infinity. The critical values, i.e.\ the images of the critical points, of $G_{a_0,b}$ are given by
$$v_{\pm}=G_{a_0,b}(c_{\pm})=bc_{\pm}^3\frac{c_{\pm}-a_0}{1-c_{\pm}a_0},$$
where the critical points $c_{\pm}$ do not depend on $b$. Moreover, if $a\neq\pm1$, at least one of the critical values, say $v$, is different from zero. Indeed, the rational maps $G_{a,b}$ have a unique preimage of zero at $z=a$. The claim holds since the critical points collide only if $a=\pm 2$ and $c_+(\pm2)=c_-(\pm2)=\pm1$. If $|b|$ is large enough  the critical value $v$ satisfies $|v|>\lambda(|a_0|+1)$ with $\lambda> 1$. As in the proof of Lemma~\ref{boundednonescapingGab1}, we conclude that the orbit of   $v$ accumulates on infinity and, therefore, for $|b|$ large enough the parameter $(a_0,b)$ is escaping.

\endproof
\section{Tongues of the Blaschke family}\label{stongues}

In this section we introduce the concept of tongue for the Blaschcke family $B_a$ restricted to the set of parameters $a$ such that $|a|\geq2$. The definition generalizes to more general almost bicritical families $f_a$  which leave the unit circle invariant and such that $f_a|_{\cercle}$ is strictly increasing of degree $2$ (c.f.\ \cite{MiRo1}). Recall that, by almost bicritical families we mean that the maps $f_a$ have at most two free critical points whilst the other ones are permanent superattracting fixed points.

Let $\{f_a\}_{a\in\Delta}$, where  $\Delta\subset\wcom $ or $\Delta\subset \real^2$, be a family of orientation preserving homeomorphisms of the unit circle. We can assign a rotation number to each of its members as follows. If $F_a$ is a lift of $f_a$, then the rotation number of $f_a$ is defined as
$$
\rho(f_a)=\lim_{n\rightarrow\infty}\frac{F_{a}^n(x)}{n},
$$

\noindent where $x\in\real$. It is well known that the limit exists and it is independent of $x$ (see e.g.\ \cite{deMVa}). Since any two lifts differ in an integer constant, this limit provides a well defined rotation number in $\real/\mathbb{Z}$. If the maps $f_a$ are  degree $2$ covers of $\cercle$, then the generalized limit does not only depend on the point $x$, but provides a semiconjugacy between $f_a$ and the doubling map. Indeed, the following lemma holds.

\begin{lemma}[{\cite[Lemmas~3.1 and 3.3]{MiRo1}}]\label{semiconj}
Let $F_a:\real\rightarrow\real$ be a continuous and increasing map depending continuously on $a$. Suppose that $F_a(x+k)=F_a(x)+2k$ for any integer $k$ and for all $x\in\real$. Then, the limit

$$
H_a(x)=\lim_{n\rightarrow\infty}\frac{F_a^n(x)}{2^n}
$$
exists uniformly on $x$.  The map $H_a$ is  increasing, continuous, depends continuously on $a$ and satisfies $H_a(x+k)=H_a(x)+k$ for any integer $k$ and for all $x\in\real$. Moreover, $H_a$ semiconjugates $F_a$ with the multiplication by $2$, i.e.\ $H_a(F_a(x))=2H_a(x)$ for any real $x$. Furthermore, if $F_a$ is increasing with respect to $a$ (for any fixed $x$), then $H_a$ is also increasing with respect to $a$.

\end{lemma}

If $F$ is a lift of a  degree 2 orientation preserving map $f:\cercle\rightarrow \cercle$, then $F$ satisfies the conditions of Lemma \ref{semiconj} and is semiconjugate to the doubling map $\theta\rightarrow 2\theta$ via  $H:\real\rightarrow\real$. Let $h$ be the  degree $1$ function of the circle which has $H$  as its lift.  Then, $h$ semiconjugates $f$ with the doubling map of the circle  $\theta\rightarrow 2\theta\;(\modul \;1)$ (equivalently given by $R_2(z)= z^2$, $|z|=1$). The following lemmas tell us that the semiconjugacy is unique  and sends periodic points to periodic points of the same period .

\begin{lemma}[{\cite{Boyl2}}]\label{unique}
Let $f$ be a degree 2 orientation preserving map of the circle. Then there exists a unique degree $1$ map $h$ of the circle which semiconjugates $f$ with the doubling map $R_2$. 
\end{lemma}

\begin{lemma}[{\cite[Lem.~3.2]{MiRo1}}]\label{periodicaperiodic}
The semiconjugating map $h$ sends points of period $k$ to points of period $k$.

\end{lemma}

Using this unique semiconjugacy $h_a$ of $B_a|_{\cercle}$, tongues can  be defined as sets of parameters for which $B_a$ has an attracting cycle in $\cercle$. Before defining them we introduce a technical lemma  and some notation. We denote by $\langle x_0\rangle=\{x_0, x_1=B_a(x_0),\cdots, x_{p-1}=B_a^{p-1}(x_0)\}$ a  period-$p$ attracting or parabolic cycle and by $A^*(\langle x_0\rangle)$ its immediate basin of attraction.

\begin{lemma}\label{one attractingpar}
A Blaschke product $B_a$, where $a\geq2$, can have at most one attracting or parabolic cycle  $\langle x_0\rangle$ in the unit circle, which has a real multiplier. If the cycle is attracting, the two critical points lie in the same connected component of $A^*(\langle x_0\rangle)$. If the cycle is parabolic, then it has multiplier $\lambda=1$ and either every point $x_n$ of the cycle lies in the boundary of a unique connected component of $A^*(\langle x_0\rangle)$ which intersects the unit circle (see Figure~\ref{bifneartip} (b)) and the two critical points lie in the same connected component of $A^*(\langle x_0\rangle)$ or there are two such components which are symmetric and which do not intersect the unit circle (see Figure \ref{bifneartip} (a)).
\end{lemma}
\proof
The multiplier is real since  $\cercle$ is invariant under $B_a$.
 It follows from the symmetry of the Blaschke family that if the two free critical points do not lie in the unit circle, then  their orbits are symmetric with respect to $\cercle$.  Hence, if one of the critical orbits accumulates on an attracting or parabolic cycle $\langle x_0\rangle$ in $\cercle$,  so does the other one. Given that any attracting or parabolic  cycle has a critical point in its immediate basin of attraction, this proves that there can be at most one attracting or parabolic cycle in the unit circle. Moreover, every connected component of $A^*(\langle x_0\rangle)$ which intersects $\cercle$  is symmetric with respect to it. If the cycle is attracting, then all connected components of $A^*(\langle x_0\rangle)$ intersect the unit circle and at least one of them contains both free critical points. This finishes the attracting case.
 
Now suppose that $\langle x_0\rangle$ is a parabolic cycle. Then it has multiplier $\lambda=1$. Indeed, since $\lambda$ is real, $\lambda=\pm 1$, but  it is positive since $B_a|_{\cercle}$ is increasing (see Lemma~\ref{increasing}). Hence, all the connected components of $A^*(\langle x_0\rangle)$ have the same period $p$ as has the cycle $\langle x_0\rangle$. If a connected component of $A^*(\langle x_0\rangle)$ intersects the unit circle, then we have a cycle of such components intersecting $\cercle$. As in the attracting case, we conclude that there is a connected component of the cycle which contains both free critical points and there can be no other parabolic cycle. Finally, if no connected component of $A^*(\langle x_0\rangle)$  intersects $\cercle$, by symmetry, every point $x_n$ of the cycle has at least $2$ such components attached to it. Since every cycle of parabolic attracting basins contains at least a critical point and there are only two free critical points, we conclude that every point of the cycle has exactly 2 components of $A^*(\langle x_0\rangle)$ attached to it.

\endproof

\begin{defi}
Assume that $B_a$, where $|a|\geq2$, has an attracting cycle $\langle x_0\rangle$ in the unit circle. The point $x_j\in\langle x_0\rangle$ such that the critical points lie in $A^*(x_j)$ is called the \emph{marked point} of the cycle. 
\end{defi}

 We generally rename the cycle so that $x_0$ denotes the marked point. Now we can formalize the concept of tongue for degree two covers of the circle. Let $H_a$ be the continuous map given by Lemma~\ref{semiconj}  which semiconjugates the lift $F_{a}$ of $B_a|_{\cercle}$ to the doubling map. Then tongues for the Blaschke products $B_a$ are defined as follows.

\begin{defin}\label{deftongue}
Let $B_a$ be a Blaschke product. We say that a parameter $a$, $|a|\geq2$, is of  \textit{type} $\tau$ if $B_a|_{\cercle}$ has an attracting cycle $\langle x_0\rangle$ and $H_a(x_0)=\tau$, where $x_0$ is the marked point point of the cycle. The \textit{tongue} $T_{\tau}$ is defined as the set of parameters $a$, $|a|\geq2$, such that $a$ is of type $\tau$.
\end{defin}

The type $\tau(a)$ is a well defined number of $\real/\mathbb{Z}$ by Lemma~\ref{unique}. Hence, we may assume that $\tau(a)\in[0,1)$. Notice that in the previous definition we use an abuse of notation on the definition, naming $x_0$ both the point in the unit circle and its lifted equivalent in the real line.  

It follows from  Lemma \ref{one attractingpar} that the tongues are disjoint. Indeed, if two different tongues would intersect, we would have parameters with two different attracting cycles in $\cercle$, which is not possible.

Given that $H_a$ sends periodic points to periodic points (see Lemma~\ref{periodicaperiodic}), any realizable type $\tau\in\cercle$ is a periodic point of the doubling map. It also follows from this and the continuity of $H_a$ with respect to parameters that tongues are open subsets of $\wcom\setminus\dis_2=\{a\in\wcom\;|\;|a|\geq2\}$. Therefore, a parameter $a$ with $|a|=2$ such that $B_a$ has an attracting cycle in $\cercle$ is said to be  in the interior of a tongue (c.f.\ Section~\ref{exttongue}).

\begin{figure}[hbt]
    \centering
    \subfigure[\scriptsize{Tongues}  ]{
    \fbox{\def\svgwidth{200pt}\subimport{/}{figure3a.pdf_tex}}}
    \hspace{0.1in}
    \subfigure[\scriptsize{Zoom in the tongues}  ]{
    \fbox{\includegraphics[width=200pt]{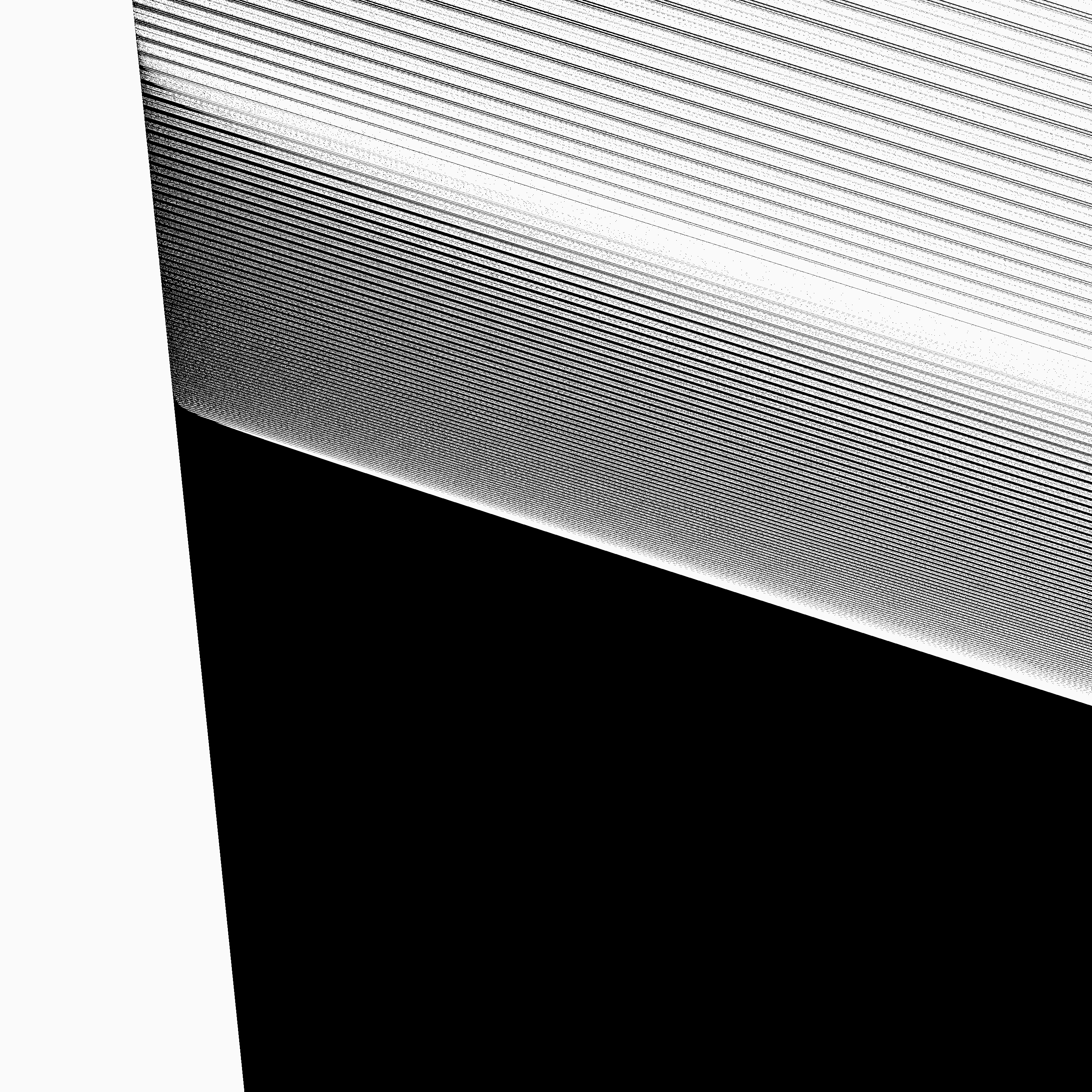}}}
  \caption{\small{In figure (a) we show the tongues of the Blaschke family for $a=re^{2\pi i \alpha}$ such that $0<\alpha<1/6$. Notice that we know, from the symmetries explained in Lemma \ref{symm}, that these parameters give complete information about the family. In figure (b) we zoom near the boundary of $T_0$. We can see how smaller tongues accumulate on it. }}   
  \label{tongues}
\end{figure}

The remaining part of the section is devoted to prove Theorem~A, which explains the basic topological properties of the tongues of the Blaschke family. Its proof splits in the next two subsections.

\subsection{Connectivity of the tongues: proof of statements (a) and (b) of Theorem~A}\label{contongues}

In this subsection we  prove  statements (a) and (b) of Theorem~A. The proof uses quasiconformal tools. We refer to \cite{Ah} and \cite{BF} for an introduction to the topic.

The proof  is inspired by Dezotti~\cite{De} and consists in performing a continuous change of the multiplier of the attracting cycle $\langle x_0\rangle\subset \cercle$.  Given a parameter $a$ of type $\tau$, with multiplier $\lambda\neq 0$, we make a quasiconformal modification of the function that  changes the multiplier to $\rho\in (0,1)$ while leaving the rest of the dynamics unchanged, obtaining a new parameter $a(\rho)$. With this modification we shall obtain a path $\gamma\subset T_{\tau}$ landing on a parameter $r_{\tau}\in T_{\tau}$ having a supperattracting fixed point. It follows from the construction that every connected component of a tongue $T_{\tau}$ contains a superatracting parameter. The proof is  finished by showing that, for every periodic point $\tau$ of the doubling map, there exists a unique parameter $r_{\tau}$ of the form $2e^{2\pi \alpha_{\tau}}$ with $\alpha_{\tau}\in[0,1/3)$ such that $B_{r_{\tau}}$ has a superatracting periodic point and has type $\tau$. Moreover,  this also shows that any tongue $T_{\tau}$ is not empty and has a unique connected component modulo the symmetries given by the third root of unity. 

We begin by changing the multiplier. The main steps of this quasiconformal construction are the following.

\begin{enumerate}[(a)]
\item First we consider a linearising map $\phi$ of $B_a^{p}$ around the attracting periodic point $x_0$, which has period $p$ and multiplier $\lambda\in(0,1)$.
\item Then we define a quasiconformal conjugacy $\mathcal{X}$ between the maps $z\rightarrow \lambda z$ and $z\rightarrow \rho z$, where $\rho\in(0,1)$. 
\item We continue by defining a $B_a-$invariant Beltrami form $\mu_{\rho}$. Around $x_0$, it is defined by pulling back the standard Beltrami coefficient $\mu_0\equiv 0$ by the quasiconformal homeomorphism $\mathcal{X}\circ\phi$. Then we spread it by the dynamics of $B_a$.
\item Finally we consider the map $\varphi_{\rho}\circ B_a\circ \varphi_{\rho}^{-1}$, where $\varphi_{\rho}$ is the integrating map of $\mu_{\rho}$ which fixes zero, infinity and $x_0$ given by the Measurable Riemann Mapping Theorem (see \cite[p. 57]{Ah}, \cite[Theorem 1.28]{BF}). This map is holomorphic and linearly conjugate to a member $B_{a(\rho)}$ of the Blaschke family.
\end{enumerate}

We now proceed to make the construction precise. Let $a\in T_{\tau}$. Recall from Definition~\ref{deftongue} that then $B_{a}|_{\cercle}$ has an attracting cycle  $\langle x_0\rangle=\{x_0,...,x_{p-1}\}\subset\cercle$ of multiplier $\lambda$.  We assume that $x_0$ lies in the component of the immediate basin which contains the critical points. Notice that $\lambda\in\real$ since $B_a|_{\cercle}$ is an endomorphism of the unit circle. Let  $\phi:U\rightarrow\dis$ be the K{\oe}nigs linearizer of $ B_a^ p$ around $x_0$, where $U$ denotes the maximal neighbourhood of $x_0$ in which $\phi$ can be taken to be conformal (see \cite[Theorem~8.2 and Lemma~8.5]{Mi1}), normalized as in the following lemma.
\begin{lemma}
The map $\phi:U\rightarrow \dis$ may be chosen to satisfy $\phi(\mathcal{I}(z))=\overline{\phi(z)}$, where $\inv(z)=1/\overline{z}$. Moreover, $U=\inv(U)$

\end{lemma}
\proof
The map $\phi$ sends invariant curves of $B_a|_{U}$ to invariant curves of  $z\rightarrow \lambda z$, which are straight lines going through $z=0$ since $\lambda$ is real. Hence,  we may assume that \linebreak $\phi(U\cap\cercle)\subset \real$ postcomposing $\phi$ with a rotation.
With the previous normalization, notice that the holomorphic map $\widehat{\phi}(z)=\overline{\phi(1/\overline{z})}$ coincides with $\phi$ on $U\cap\cercle$ and, therefore, it equals $\phi$. The symmetry of $U$ follows from the symmetry of $\phi$.

\endproof

In the following lemma we introduce  a quasiconformal map $\mathcal{X}$ which is used to change the multiplier of an attracting cycle. We refer to \cite{De} and \cite{BF} for a proof.

\begin{lemma}\label{deform}
Let
$$\mathcal{X}:\com^*\rightarrow \com^*$$
$$z\rightarrow |z|^{\alpha}z,$$

\noindent where $\alpha\in(-1,\infty)$ and let $0<R<1$ and $0<r<1$. Then the following hold.

\begin{enumerate}[(a)]
\item The Beltrami form $\mu_{\mathcal{X}}=\mathcal{X}^*\mu_0$, where $\mu_0$ denotes the Beltrami form of the standard complex structure and $^*$ denotes the pull back operation,  satisfies
$$\mu_{\mathcal{X}}=\frac{\partial \mathcal{X}/\partial \overline{z}}{\partial \mathcal{X}/\partial z}=\frac{\alpha}{2+\alpha}\frac{z}{\overline{z}},$$
$$||\mu_{\XX}||_{\infty}=\left|\frac{\alpha}{2+\alpha}\right|.$$
\item $\XX$ is invertible and satisfies $\XX(r e^{2\pi i\theta})= \xi(r)e^{2\pi i\theta}$, where $\xi(r)=r^{\alpha+1}$.
\item Let $\alpha=\frac{\log r}{\log R}-1$. Then $\XX$ sends the disk $\dis_R$ of radius $R$ to the disk $\dis_r$ of radius $r$ and, moreover,
$$||\mu_{\XX}||_{\infty}=\frac{|1-\log r/\log R|}{1+\log r/\log R}<1.$$
\item Let $\lambda\in (0,1)$. Let $\rho=\XX(\lambda)=\lambda^{1+\alpha}$. Then  $\XX$ conjugates the map multiplication by $\lambda$ with the map multiplication by $\rho$, i.e.\ the following diagram is commutative
\begin{equation*}\begin{CD}
\dis_R @> z\rightarrow\lambda z >> \dis_{\lambda R}\\
@V  \XX VV @VV \XX  V\\
\dis_r@>  z\rightarrow\rho z>>\dis_{\rho r}.
\end{CD}\end{equation*}

\end{enumerate}

 \end{lemma}

 Notice that, since $\lambda\in (0,1)$ and $\alpha\in(-1,\infty)$, we have that $0<\rho<1$.

\begin{lemma}
The Beltrami form given by $\mu_{\rho}=\mu_{\XX\circ\phi}=\phi^*\mu_{\XX}$ depends analytically on $\rho$ and is invariant under $B_a^{p}$. Moreover, $\mathcal{I}^*\mu_{\rho}=\mu_{\rho}$ on $U$.
\end{lemma}

\proof

The analytic dependence with respect to $\rho$ is obtained from the explicit expression of $\mu_{\mathcal{X}}$ and the fact that $\alpha=\frac{\log \rho}{\log \lambda}-1$.  Invariance under $B_a^{ p}$ follows from the next commutative diagram.

\begin{equation*}\begin{CD}
(U, \mu_{\XX\circ\phi})@> \phi >>(\dis_R, \mu_{\mathcal{X}}) @> \XX >> (\dis_r, \mu_0)\\
@V B_a^{ p} VV  @V z\rightarrow\lambda z VV @VV z\rightarrow\rho z V\\
(B_a^{ p}(U), \mu_{\XX\circ\phi})@> \phi >>(\dis_{\lambda R}, \mu_{\mathcal{X}})@>  \XX >>(\dis_{\rho r}, \mu_0).
\end{CD}\end{equation*}

To see that $\mathcal{I}^*\mu_{\rho}=\mu_{\rho}$, we have to check that $\mu_{\rho}(z)=\overline{\mu_{\rho}(1/\overline{z})}z^2/\overline{z}^2$ (see \cite[Definition 1.3]{BF}). From the explicit expression of the Beltrami form that we get from the pullback (see \cite[Equation (1.9)]{BF}), we have
$$\mu_{\rho}(z)=\frac{\alpha/2}{1+\alpha/2}\frac{\phi(z)}{\overline{\phi(z)}}\cdot \frac{\overline{\phi'(z)}}{\phi'(z)}.$$

The result follows since $\phi(1/\bar{z})=\overline{\phi(z)}$ and $\phi'(z)=-\overline{\phi'(1/\bar{z})}/z^2$.

\endproof

Once we have this Beltrami form given by $\mu_{\rho}$ in $U$, we spread it to $\wcom$ by defining:

$$\mu_{\rho}=\left\{\begin{array}{lcl}
\mu_{\rho} &  \mbox{ on } & U\\
(B_a^n)^*\mu_{\rho} &\mbox{ on } &  B_a^{-n}(U)\setminus B_a^{-n+1}(U),\mbox{ for } n>1\\
\mu_0 &   \mbox{otherwise,} &
\end{array}\right. $$

\noindent where $\mu_0\equiv0$. Then $\mu_{\rho}$ is well defined since it is $B_a^{ p}$-invariant. It depends analytically on $\rho$ since we are pulling back by a holomorphic map an almost complex structure which  depends analytically on $\rho$. Furthermore, since $\mu_{\rho}|_{U}$ is symmetric with respect to $\cercle$ and $\mu_{\rho}$ is defined recursively by pulling back by $B_a$, which is also symmetric with respect to $\cercle$, $\mu_{\rho}$ also satisfies $\mathcal{I}^*\mu_{\rho}(z)=\mu_{\rho}(z)$ for all $z\in\com$ (see \cite[Exercise 1.2.5]{BF}).

Since $\mu_{\rho}$ is built pulling back $\mu_{\mathcal{X}}$ by holomorphic mappings and we have that $||\mu_{\mathcal{X}}||_{\infty}<1$,  we also have that $||\mu_{\rho}||_{\infty}<1$. Let $\varphi_{\rho}$ be the integrating map obtained by applying the Measurable Riemann Mapping Theorem (see \cite[p. 57]{Ah}, \cite[Theorem 1.28]{BF}) such that it fixes $0$, $x_0$ and $\infty$. The next lemma follows from uniqueness of the integrating map

\begin{lemma}\label{simmetry}
The integrating map $\varphi_{\rho}$ is symmetric with respect to $\cercle$.

\end{lemma}
\proof
Consider the map $\widetilde{\varphi}_{\rho}(z)=\inv\circ\varphi_{\rho}\circ\inv(z)$, where $\inv(z)=1/\overline{z}$. By construction it fixes $0$, $x_0$ and $\infty$. Moreover, by symmetry of $\mu$ with respect to the unit circle we have $\widetilde{\varphi}_{\rho}^*\mu=\mu$. Then, by uniqueness of the integrating map we conclude that $\widetilde{\varphi}_{\rho}(z)=\varphi_{\rho}(z)$.
\endproof

Once we have $\varphi_{\rho}$, we can build our new rational map.

\begin{propo}
The map $\varphi_{\rho}\circ B_a \circ \varphi_{\rho}^{-1}$ is  a rational map of degree $4$ of the form $B_{a_{\rho},t_{\rho}}$ as in (\ref{blasformula2}), where $t_{\rho}\in\real$ and $a_{\rho}\in\com$. The parameters $a_{\rho}$ and $t_{\rho}$ depend continuously on $\rho$. Moreover,  the attracting fixed point $x_0$ of $B_{a_{\rho},t_{\rho}}^p$ has multiplier $\rho$. 
\end{propo}
\proof

By construction,  the quasiregular map $\varphi_{\rho}\circ B_a \circ \varphi_{\rho}^{-1}$ preserves the standard complex structure. Consequently, it is a holomorphic map of the Riemann sphere onto itself by Weyl's Lemma (see \cite[p.\ 16]{Ah}, c.f.\ \cite[Theorem 1.14]{BF}). It is of the form $B_{a_{\rho},t_{\rho}}$ since it has local degree $3$ around $0$ and $\infty$, it has global degree $4$ and it is symmetric with respect to the unit circle (see \cite[Lem.~15.5]{Mi1}).

We now check the dependence on parameters. Recall that $B_a$ has a unique critical point $c_+\in\com\setminus\dis$ (see Section \ref{introblas0}). Since $\mu_{\rho}$ depends real analytically on $\rho$, the integrating map $\varphi_{\rho}(z)$ depends real analytically on $\rho$ for all $z\in \com$ by the analytic dependence on parameters of the Measurable Riemann Mapping Theorem. Therefore, the critical point $\varphi_{\rho}(c_+)\in\com\setminus\dis$ of $B_{a_{\rho},t_{\rho}}$ depends real analytically on $\rho$. We conclude from Lemma \ref{continuitya} that  $a_{\rho}$  depends continuously on $\rho$.  The parameter $t_{\rho}$ is continuously determined by the parameter $a_{\rho}$ together with the image of a point $z_0\neq0,\infty$. Given that the map $B_{a_{\rho},t_{\rho}}(x_0)=\varphi_{\rho}\circ B_a(x_0)$ depends real analytically on $\rho$, we conclude that $t_{\rho}$ also depends continuously on $\rho$.

Finally, let $\widetilde{\phi}(z)=\mathcal{X}\circ \phi\circ \varphi_{\rho}^{-1}(z)$, where $z\in\varphi_{\rho}(U)$. By construction $\widetilde{\phi}$ is a quasiconformal map which conjugates $B_{a_{\rho},t_{\rho}}^p$ around $x_0$ to the map $z\rightarrow \rho z$. Since it preserves the standard complex structure, it is a conformal map by Weyl's Lemma. Hence, it is the linearizing function and $\rho$ is the multiplier of the new cycle.

\endproof

We know from Lemma \ref{conjblas} that $B_{a_{\rho}, t_{\rho}}$ is conjugate to $B_{a_{\rho} e^{-i t_{\rho}/3},0}=B_{a_{\rho} e^{-i t_{\rho}/3}}$ by a conjugacy $\mathcal{L}_{\rho}$. This gives us a continuous curve in the set of parameters of our Blaschke family. Indeed, we have a (real analytic) curve  $\gamma:(0,1)\rightarrow \com\setminus \dis_2$ defined as $\gamma(\rho)=a_{\rho} e^{-i t_{\rho}/3}=a(\rho)$.

\begin{lemma}\label{type}
If the parameter $a$ has type $\tau$, then, for all $\rho\in(0,1)$, the parameter $a(\rho)$ has type $\tau$. 
\end{lemma}
\proof

We have already seen that $B_{a(\rho)}$ has an attracting cycle, so it has type $\tau'$. It is easy to check that the map 
$\varrho(z)=\mathcal{L}_{\rho}\circ \varphi_{\rho}(z)$
conjugates $B_{a(\rho)}$ and $B_{a}$, sending the marked periodic point $x_0$ to $\varrho(x_0)$. Moreover, the  continuous map $H_a\circ\varrho^{-1}$  semiconjugates $B_{a(\rho)}$ with the doubling map. Therefore, we have that

$$\tau'=H_a\circ\varrho^{-1}(\varrho(x_0))=H_a(x_0)=\tau.$$
\endproof

Now we have a path $\gamma(\rho)=a(\rho)$ defined for $\rho\in(0,1)$ which gives, for each $\rho$, a parameter $a(\rho)\in T_{\tau}$. We want to prove that this path lands at a single point when $\rho\rightarrow 0$. Note that $|a(\rho)|\rightarrow 2$ when $\rho\rightarrow 0$ since $|B_{a}'|_{\cercle}|>C>0$ when $|a|>2+\epsilon$, where $\epsilon>0$ and $C=C(\epsilon)$ is a constant. It follows from the continuous dependence on $a$ of the semiconjugacy $H_a$ (see Lemma \ref{semiconj}) that any limit point of $\gamma$  has a superattracting fixed point of period $p$.

Let $\omega$ be the limit set of $\gamma(\rho)$ when $\rho\rightarrow 0$. Since $\omega=\bigcap_n\overline{\gamma((0,1/n))}$ is a decreasing intersection of connected compact sets, we conclude that it is a connected set of parameters $a$ such that $|a|=2$.

We restrict now to parameters $a$ such that $|a|=2$.  Let $a=2e^{2\pi i \alpha}$. Throughout the rest of the proof it will be convenient to work with $B_{2,3\alpha}$ (Equation (\ref{blasformula2})) as in Section \ref{prelimparametrization} so as to use Lemma~\ref{increasing}. This map is conformally conjugate to $B_a=B_{a,0}$  (Equation (\ref{blasformula})) by the rotation $\mathcal{L}(z)=e^{-2\pi i \alpha}z$ (see Lemma~\ref{conjblas}). For $\alpha\in[0,1/3)$ we have a one to one correspondence between the parameters $2e^{2\pi i \alpha}$ and the parameters $(2, 3\alpha\;(\modul \;1))$. Notice that, for all $\alpha$, $B_{2,3\alpha}$ has a unique critical point at $c=1$. Indeed, for $a=2e^{2\pi i \alpha}$, the two critical points collapse in $c(\alpha)=e^{2\pi i \alpha}$, which is sent to $c=1$ by the conjugacy $\mathcal{L}(z)$. 

Assume that $\omega$ is not a single parameter. Then, we have  an interval of parameters with a superattracting periodic cycle. Therefore, the critical point $c=1$ is periodic in this interval of parameters,  i.e.\ $B_{2,3\alpha}^{ p}(1)=1$  for all parameters $ (2,3\alpha)\in\omega$. This is impossible since $B_{2,\alpha}$  is strictly increasing with respect to $\alpha$ (see Lemma \ref{increasing}). Hence, $\omega$ is a single parameter $r_{\tau}=\lim_{\rho\rightarrow 0}a(\rho)$. Notice also that this parameter has type $\tau$. Indeed,  given the fact that it has a superattracting periodic point, it belongs to a tongue of type $\tau'$. Since tongues are open sets in $\com\setminus\dis_2=\{a\; | \; |a|\geq2\}$, we conclude that any curve of parameters contained in $\com\setminus\dis_2$ and landing in $\omega$ necessarily intersects $T_{\tau'}$. We conclude that $\tau'=\tau$ since $a(\rho)$ have type $\tau$ for all $\rho\in(0,1)$ by Lemma \ref{type}.

In order to finish the proof of statements (a) and (b) of Theorem~A we have to show that the limit does not depend on the initial parameter $a\in T_{\tau}$.  We use the following lemma.

\begin{lemma}\label{increasing2}
Let $g_{\alpha}(x):=B_{2,3\alpha}|_{\cercle}(x)$, $x\in\cercle$. Then, for any $p\in \mathbb{N}$, the mapping \linebreak $\alpha\rightarrow g_{\alpha}^{p}(1)\in \cercle$, $\alpha\in[0,1/3)$, is strictly increasing and of degree $2^p -1$.

\end{lemma}
\proof
The map $g_{\alpha}$ is increassing by Lemma \ref{increasing}. We only have to prove that $g_{\alpha}$ has degree $2^p -1$. The lift of $g_{\alpha}$ is given by

 $$h_{\alpha}(x)=3x+3\alpha+\frac{1}{2\pi i}\log\left(\frac{e^{2\pi i x}-2}{1-2e^{2\pi i x}}\right).$$

The result is true for $p=1$ in the sense that $h_{\alpha+1/3}(x)=h_{\alpha}(x)+1$. By induction over $p$ and using that $h_{\alpha}^p(x+1)=h^p_{\alpha}(x)+2^p$ ($g_\alpha$ has degre two as a circle map) we have

$$h_{\alpha+1/3}^{(p+1)}(x)=h_{\alpha+1/3}^{p}(h_{\alpha+1/3}(x))=h_{\alpha+1/3}^{ p}(h_{\alpha}(x)+1)=h_{\alpha+1/3}^{ p}(h_{\alpha}(x))+2^p$$
$$=h_{\alpha}^{ p}(h_{\alpha}(x))+2^p-1+2^p=h_{\alpha}^{ p+1}(x)+2^{p+1}-1.$$
\endproof

 It follows from this lemma that $g_{\alpha}^p$ has exactly $2^p-1$ parameters $\alpha$,  $\alpha\in[0,1/3)$, such that the critical point is periodic of period dividing $p$. Indeed, for every natural $k\in\{0,1,...,2^p-2\}$, there exists a unique $\alpha_{p,k}\in[0,1/3)$ such that $h_{\alpha_{p,k}}^p(0)=0+k$. It can be computed using the expression of the semiconjugacy $H_a$ (see Lemma \ref{semiconj}) that a parameter $a_{k,p}=2e^{2\pi i\alpha_{k,p}}$ has type $\tau(a_{k,p})=k/(2^p-1)$ (c.f.\ \cite[Lem.~4.2]{De}). Since the expression in the form $k/(2^p-1)$  of a periodic point $\tau$ of period $p$ of the doubling map is unique, we conclude that, for a fixed a type $\tau$, there exists a unique parameter  $\alpha_{p,k}\in[0,1/3)$ which has a superattracting cycle of type $\tau$. Hence, we can also conclude that no tongue $T_{\tau}$ is empty. It also follows from this that the limit $\omega$ is unique up to conjugacy since, for a fixed a type $\tau$, there exists a unique possible limit. This finishes the proof of statements (a) and (b) of Theorem~A.

\begin{rem}
We have not proven that the three connected components of a tongue $T_{\tau}$ are actually different. A priori a tongue $T_{\tau}$ could have a single connected component that would contain the three symmetric roots. However, this is not possible due to the description of the boundaries of the tongues given by statement (d) of Theorem~A (c.f.\ Theorem~\ref{simplygraph}), whose proof is independent of the one of statements (a) and (b).
\end{rem}

\subsection{Boundary of the tongues: proof of statements (c) and (d) of Theorem~A}\label{Boundary Tongues}

 The goal of this subsection is to prove statements (c) and (d) of Theorem~A (c.f.\ \cite{MiRo1,MiRo2}). In order to do so we describe the boundary of the tongues and the parameters therein. 
 The following is a direct consequence of Lemma~\ref{increasing}, which states that  there cannot be any attracting periodic cycle in $\cercle$ if $|a|\geq3$. 
 
\begin{propo}\label{boundedtongues}
If $a\in T_{\tau}$ then $|a|<3$. Consequently, tongues are bounded.
\end{propo}

The following result is a corollary of Lemma~\ref{one attractingpar},  which states that a parabolic cycle of $B_{a}|_{\cercle}$ can have at most one periodic petal intersecting $\cercle$.
\begin{co}\label{paroneside}
A parabolic cycle $\langle x_0\rangle\in\cercle$ of $B_a$ cannot be attracting from both sides in $\cercle$. More precisely, if $x_0$ is a parabolic fixed point of $B_a^p$, there cannot exist a neighborhood $U$ of $x_0$ in $\cercle$ such that all  the points  $U$ are attracted to $x_0$ under iterations of $B_a^p$.
\end{co} 
 
The next lemma states that parameters on the boundary of tongues correspond to maps having parabolic cycles. 
 
\begin{lemma}\label{boundaryparabolic}
If $a$ belongs to the boundary of a tongue, then $B_a$  has a parabolic cycle of multiplier $1$.
\end{lemma}
\proof
Let $a_0\in\partial T_{\tau}$. Then, there exists a sequence of parameters $a_n\in T_{\tau}$, $n\in \nat$, such that $a_n$ accumulate on $a_0$. Let $x_n$ be the attracting periodic point of $B_{a_n}$ having the critical points in its immediate basin of attraction. Since $\cercle$ is compact, we may assume that $x_n$ converge to a point $x_0\in\cercle$. Since $B_a$ depends continuously on $a$, we conclude that $x_0$ is a periodic point of $B_{a_0}$.  The multiplier of $x_0$ has to be  $1$. Indeed, it is real since $x_0\in\cercle$, it is positive since $B_a|_{\cercle}$ is increasing and it cannot be smaller than $1$ because otherwise would belong to the interior of a tongue $T_{\tau'}$, which is impossible since tongues are disjoint.

\endproof

From now on it will be convenient to work with the alternative parametrization of the Blaschke family  $g_{r,\alpha}:=B_{r,3\alpha}|_{\cercle}$ as in Section~\ref{prelimparametrization}. We consider the parameter space $(r,\alpha)$ with $r\geq 2$ and $\alpha\in \real/\frac{1}{3}\mathbb{Z}$ instead of $a\in\com$ with $|a|\geq2$. The main reason to use this alternative parametrization is to use the monotonicity with respect to $\alpha$ given by Lemma~\ref{increasing}. 

\begin{defin}
We say that a parameter $(r_0,\alpha_0)$, $r_0>2$, in a boundary of a tongue $T_{\tau}$ of period $p$ belongs to the \textit{left boundary} of the tongue if there exists an $\epsilon >0$ such that, for all $0<\alpha<\epsilon$, $(r_0,\alpha_0+\alpha)$ belongs to the tongue  and $(r_0,\alpha_0-\alpha)$ does not belong to it. Conversely, we say that it belongs to the \textit{right boundary} if there exists an $\epsilon >0$ such that, for all $0<\alpha<\epsilon$, $(r_0,\alpha_0-\alpha)$ belongs to the tongue  and $(r_0,\alpha_0+\alpha)$ does not belong to it. Finally, we say that it belongs to a \textit{tip of the tongue} if there exists an $\epsilon>0$ such that, for all $\alpha\in (-\epsilon,0)\cup(0,\epsilon)$,  $(r_0,\alpha_0+\alpha)$ does not belong to the tongue.
\end{defin}

Using Lemma~\ref{increasing}, we have the following result (c.f. \cite[Lem.~4.1]{MiRo1}).

\begin{figure}[hbt]
    \centering
    \subfigure[\scriptsize{Case (a)} ]{
     \put(0,0){\includegraphics[width=\unitlength]{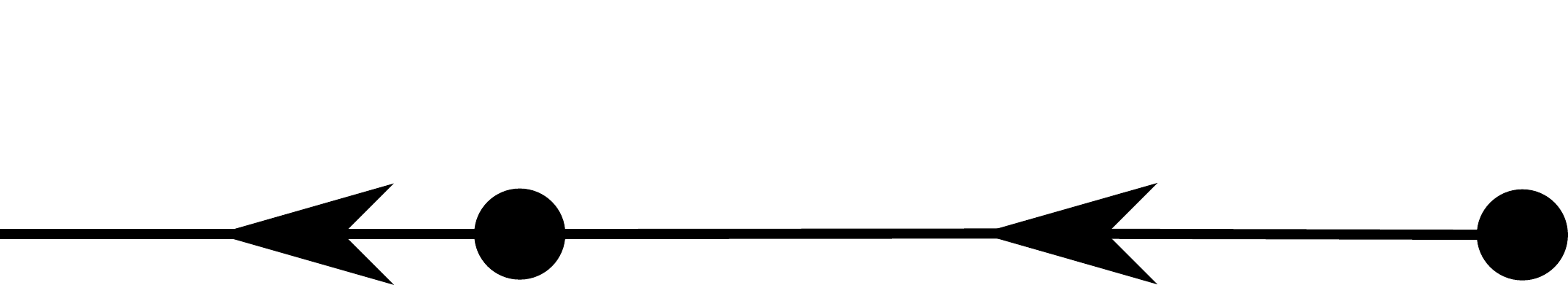}}%
    \put(60.694198,10.1374062){\color[rgb]{0,0,0}\makebox(0,0)[lb]{\smash{$J$}}}%
  \includegraphics[width=90pt]{figure4a.pdf}}
    \hspace{0.1in}
    \subfigure[ \scriptsize{Case (b)} ]{
      \put(0,0){\includegraphics[width=\unitlength]{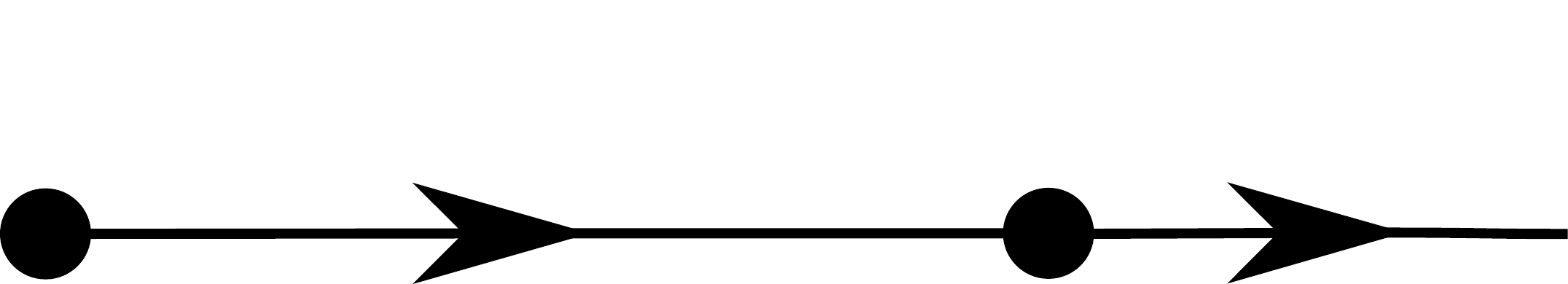}}%
    \put(25.694198,10.1374062){\color[rgb]{0,0,0}\makebox(0,0)[lb]{\smash{$J$}}}%
    \includegraphics[width=90pt]{figure4b.pdf}}
    \hspace{0.1in}  
 \subfigure[\scriptsize{Case (c)} ]{
      \put(0,0){\includegraphics[width=\unitlength]{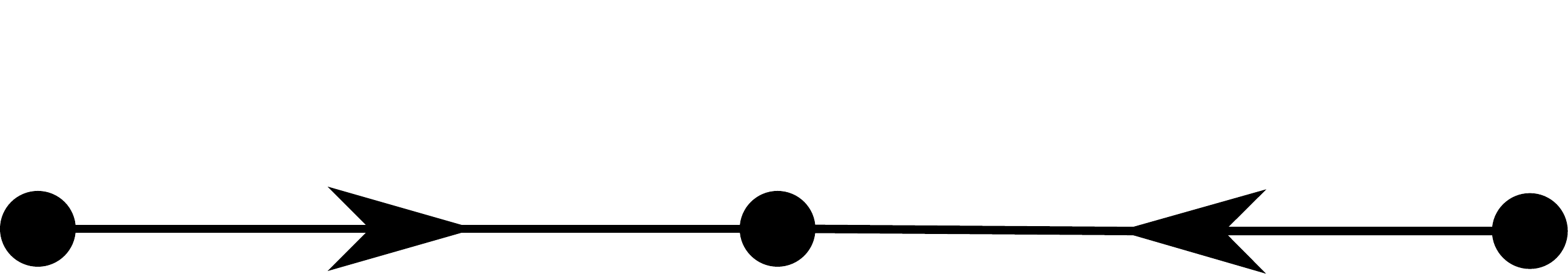}}%
    \put(40.694198,10.1374062){\color[rgb]{0,0,0}\makebox(0,0)[lb]{\smash{$J$}}}%
     \includegraphics[width=90pt]{figure4c.pdf}}
    \hspace{0.1in}  
   \subfigure[\scriptsize{Case (d)}]{
     \put(0,0){\includegraphics[width=\unitlength]{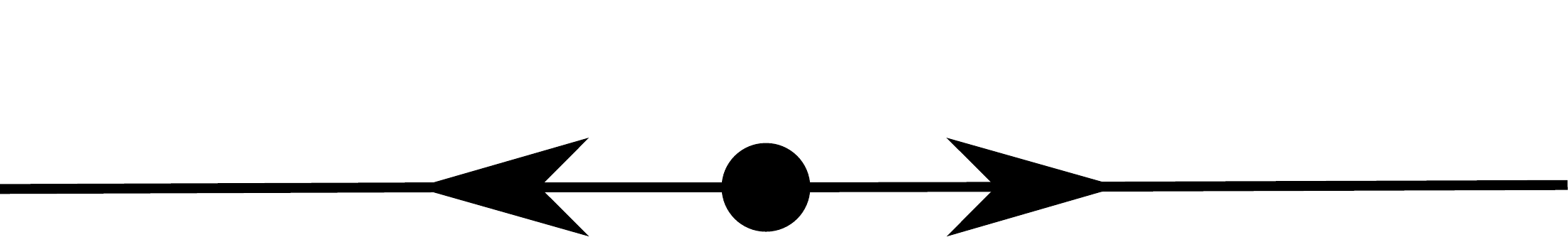}}%
    \put(29.694198,10.1374062){\color[rgb]{0,0,0}\makebox(0,0)[lb]{\smash{$J=x$}}}%
    \includegraphics[width=90pt]{figure4d.pdf}}
    \caption{\small{The four different behaviors which may occur in Lemma \ref{clasification0}.}} 
    \label{esquema}
\end{figure}

\begin{lemma}\label{clasification0}
Let $x_0$ be an attracting or parabolic periodic point of $g_{r_0,\alpha_0}$ of period $p$ and let $H_{r_0,\alpha_0}$ be the semiconjugacy between $g_{r_0,\alpha_0}$ and the doubling map given by Lemma~\ref{semiconj}. Let $J$ be the set of points $x\in\cercle$ which are sent by $H_{r_0,\alpha_0}$ to the same point as $x_0$, i.e.\ $J=\{x\;|\;H_{r_0,\alpha_0}(x)=H_{r_0,\alpha_0}(x_0)\}$. Then,  either $J$ is a connected closed interval or it consists of a single point. Moreover, $g_{r_0,\alpha_0}|_J$ is a homeomorphism, the endpoints of $J$ are fixed points of $g_{r_0,\alpha_0}^p$, and one of the following holds (see Figure \ref{esquema}).

\begin{enumerate}[(a)]
\item $J$ is an interval. The left endpoint of $J$ is parabolic, topologically attracting from the right  and repelling from the left, the right endpoint is repelling and there are no other fixed points of $g_{r_0,\alpha_0}^p$ in $J$. In this case $(r_0,\alpha_0)$ belongs to the left boundary of the tongue.
\item $J$ is an interval. The right endpoint of $J$ is parabolic, topologically attracting from the left  and repelling from the right, the left endpoint is repelling and there are no other fixed points of $g_{r_0,\alpha_0}^p$ in $J$. In this case $(r_0,\alpha_0)$ belongs to the right boundary of the tongue.
\item $J$ is an interval. Both endpoints of $J$ are repelling, there is an attracting fixed point of $g_{r_0,\alpha_0}^p$ in $J$ and there are no other fixed points of $g_{r_0,\alpha_0}^p$ in $J$.  this case $(r_0,\alpha_0)$ belongs to the interior of the tongue.
\item $J$ consists of a parabolic periodic point which is topologically repelling in $\cercle$. 
\end{enumerate}

\end{lemma}

\proof
Throughout the proof we consider the points of $\cercle$ oriented anticlockwise. The fact that $J$ is either a single point or a closed interval follows since  $H_{r_0,\alpha_0}$ is an increasing continuous function (see Lemma~\ref{semiconj}). Since $g_{r_0,\alpha_0}|_{\cercle}$ is an strictly increasing function (see Lemma~\ref{increasing}), we have  that $g_{r_0,\alpha_0}|_{\cercle}$ is a local homeomorphism around each point. It follows  that $g_{r_0,\alpha_0}|_J$ is a homeomorphism and its endpoints are fixed points of $g_{r_0,\alpha_0}^p$.

It follows from Lemma~\ref{one attractingpar}  that only these four cases can occur. It states that $g_{r_0,\alpha_0}|_{\cercle}$ can have at most one attracting or parabolic cycle  and that a para\-bolic cycle cannot be topologically attracting from both sides. 

It is left to see that case (a) corresponds to the left boundary of the tongue whereas case (b) corresponds to the right boundary. We prove it for case (a). Case (b) is analogous.

Assume that $(r_0,\alpha_0)$ satisfy the hypothesis of case (a). Then, there exists a periodic point of period $p$ which is repelling from the left and attracting from the right. Recall that, from Lemma~\ref{increasing}, we have that $g_{r,\alpha}^p$ is strictly increasing with respect to $\alpha$ for any $p\in\nat$.

We first prove that there exists an $\epsilon>0$ such that if $0<\alpha<\epsilon$, then  $g_{r_0,\alpha_0+\alpha}$ has an attracting cycle of period $p$. Let $x$ be a point of the parabolic cycle and let $y$ be a point in the immediate basin of attraction of $x$. Then $x<y$  and $g_{r_0,\alpha_0}^p(y)<y$. Since $g_{r,\alpha}^p$ is strictly increasing with respect to $\alpha$, there exists an $\epsilon>0$ such that if $0<t<\epsilon$, then $g_{r_0,\alpha_0+t}^p(x)>x$ and $g_{r_0,\alpha_0+\alpha}^p(y)<y$. Hence, there has to be a topologically attracting periodic point of period $p$ between $x$ and $y$. Since, by Corollary~\ref{paroneside}, a parabolic periodic point cannot be attracting on both sides, we get that this topological attractor located between $x$ and $y$ is an attractor.

Now we have to see that there exists an $\epsilon>0$ such that if $0<\alpha<\epsilon$, then $g_{b_0,\alpha_0-\alpha}$ has  no periodic attracting cycles of period $p$. Since $x$ is repelling from the left and attracting from the right, there exists a $\delta>0$ such that if $y\in(x-\delta,x+\delta)$,  then $g_{r_0,\alpha_0}^p(y)\leq y$. Using that $g_{r_0,\alpha_0+\alpha}$ is strictly increasing with respect to $\alpha$,  we can take $\epsilon_1 > 0$ such that if $0 < \alpha < \epsilon_1$, then $g_{r_0,\alpha_0}^p(y)< y$ for all $y\in(x-\delta,x+\delta)$. Doing the same around all of the points $x_m$ of the parabolic cycle we obtain $\delta'>0$ and $\epsilon_2>0$ such that if $0<\alpha<\epsilon_2$, then $g_{r_0,\alpha_0}^p(y)< y$ for all $y\in(x_m-\delta',x_m+\delta')$. Hence, we have erased the periodic points of period $p$ in a $\delta'$-neighborhood $U$ of our cycle.  Since $g_{r_0,\alpha_0}^p$ has finitely many fixed points in $\cercle\setminus U$, all of them repelling, we can take $\epsilon_3<\epsilon_2$ such that if $0<\alpha<\epsilon_3$, then  $g_{r_0,\alpha_0-\alpha}^p$ has no attracting fixed point at all.

\endproof

We want to prove that case (d) of Lemma~\ref{clasification0} corresponds to a tip of the tongue (see Proposition~\ref{clasification2}). We first introduce some auxiliary lemmas. The following lemma corresponds to Lemma 3.1 in \cite{MiRo2}.
\begin{lemma}\label{auxiliar}
Let $U$ be a neighborhood of the origin in $\real^2$ and let $F:U\rightarrow \real$ be a real analytic function. Set $f_t(x)=F(t,x)$. Assume that $f_0$ has a topologically repelling fixed point at $x=0$ and that

$$\frac{\partial F}{\partial t}(0,0)\neq 0,$$
 Then there are open intervals $I,J$ containing $0$ such that $I\times J\subset U$ and for every $t\in I$ the map $f_t$ has exactly one fixed point $x\in J$. Moreover,  if $t\neq 0$, then the fixed point has multiplier $\lambda>1$.

\end{lemma}

We use the previous technical lemma in the following result.

\begin{lemma}\label{flat}
Consider a one parameter subfamily $f_t=g_{r(t), \alpha(t)}$ of the Blaschke family such that $r(t)$ and $\alpha(t)$ depend analytically on $t$. Assume that $f_{t_0}^p$ has a topologically repelling parabolic fixed point $x_0$ and $\frac{\partial G}{\partial t}(t_0, x_0)\neq 0$, where $G(t,x)=f_t^p(x)$. Then there exists  $\epsilon>0$ such that if $t\in(-\epsilon,0)\cup(0,\epsilon)$ then $f_{t_0+t}^p$ has no attracting or parabolic fixed point.

\end{lemma}
\proof
By Lemma~\ref{auxiliar}, there exists $\epsilon_1>0$ and a neighborhood $U_1$ of $x_0$ such that if $t-t_0\in(-\epsilon_1,0)\cup(0,\epsilon_1)$, then $f_t^p$ has no attracting or parabolic fixed point in $U_1$. Now, as in proof of Lemma~\ref{clasification0}, we can perform the same argument around the other $p-1$ points of the parabolic cycle, obtaining an $\epsilon_2>0$ and a neighborhood $U$ of $\{x_0,...,x_{p-1}\}$ such that if $t-t_0\in(-\epsilon_2,0)\cup(0,\epsilon_2)$ then $f_t^p$ has no attracting or parabolic  fixed point in $U$. Since $f_t^p$ has only finitely many fixed points in $\cercle\setminus U$, all of them repelling, we can take $\epsilon_3<\epsilon_2$ such that if $t-t_0\in(-\epsilon_3,0)\cup(0,\epsilon_3)$ then $f_t^p$ has no attracting or parabolic fixed point at all. 
\endproof

This result gives us directly the following proposition.
\begin{propo}\label{clasification2}
A parameter $(r,\alpha)$ in a boundary of a tongue for which case (d) of Lem\-ma~\ref{clasification0} occurs, is a tip of that tongue.
\end{propo}

From Lemma \ref{clasification0} and Theorem \ref{clasification2} we obtain the next corollary.
 
 \begin{co}\label{clasification3}
Any parameter $(r_0,\alpha_0)$ of the boundary of a tongue $T_{\tau}$ either belongs to the right or the left boundary of the tongue or is a tip of the tongue.
\end{co}

We now prove the remaining statements of Theorem~A. The following theorem proves statement (c).

\begin{teor}\label{simplygraph}
Given $r_0\geq 2$, the intersection of any connected component of a tongue $T_{\tau}$ with the parameter circle $|a|=r_0$ is connected. In particular, every connected component of a tongue is simply connected.
\end{teor}
\proof
Assume that the intersection of a connected component of a tongue $T_{\tau}$ with the parameter circle $|a|=r_0$ is not connected. Then, there  exists a parameter $(r_0,\alpha_0)$ and $\epsilon>0$ such that, for any $\alpha\in(-\epsilon,0)\cup(0, \epsilon)$, the parameter $(r_0,\alpha_0+\alpha)$ belongs to the tongue. That would imply that $(r_0,\alpha_0)$ is a parameter in the boundary of $T_{\tau}$ which does not belong to the right or left boundary and which is not a tip of the tongue, contradicting Corollary~\ref{clasification3}.

\endproof

Finally, we prove statement (d).

\begin{proof}[Proof of statement (d)]

It follows from  Theorem~\ref{simplygraph} that the left and right boundaries are well defined curves. Both the left and the right boundary begin in two different parameters $a_-$ and $a_+$ with $|a_-|=2=|a_+|$. Both boundaries are bounded by Proposition \ref{boundedtongues} and hence they have to end at a point where they intersect, which is a tip. The only thing which is left to see is that the boundary of a tongue cannot be flat for any $r_0$, i.e.\ we have to see that the intersection of the boundary of $T_{\tau}$ with any parameter circle $|a|=r_0$ does not contain an interval of parameters. Notice that, by Theorem~\ref{simplygraph}, neither the left nor the right boundaries can have local maximums. The points of such  interval cannot be of the left boundary or the right boundary by definition. Hence, the parameters of this open interval are tips of the tongue. Therefore, we would have an $\epsilon>0$ such that $B_{r_0,\alpha_0+\alpha}$ has a topologically repelling parabolic fixed point for all $|\alpha|<\epsilon$ and for some $\alpha_0$ and $r_0$. However, this would contradict Lemma~\ref{flat}.

\end{proof}

\section{Bifurcations around the tip of the  tongues: proof of Theorem~B}\label{tiptongue}

\begin{figure}[p]
    \centering
    \subfigure[\scriptsize{Dynamical plane of $B_{3}$.}  ]{
    \includegraphics[width=209pt]{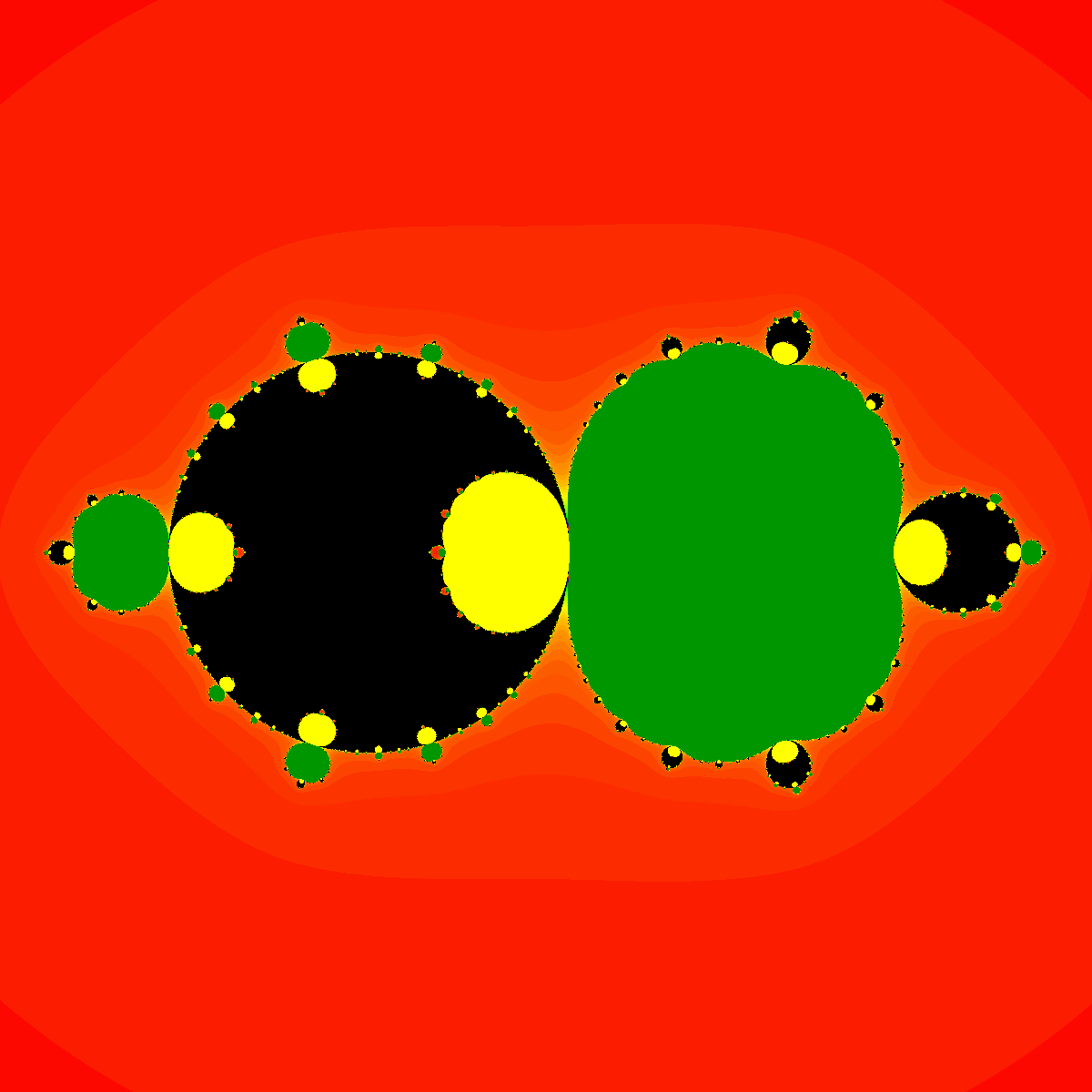}}
    \hspace{0.1in}
    \subfigure[\scriptsize{Dynamical plane of $B_{2.65675+0.0389604i}$.}  ]{
    \includegraphics[width=209pt]{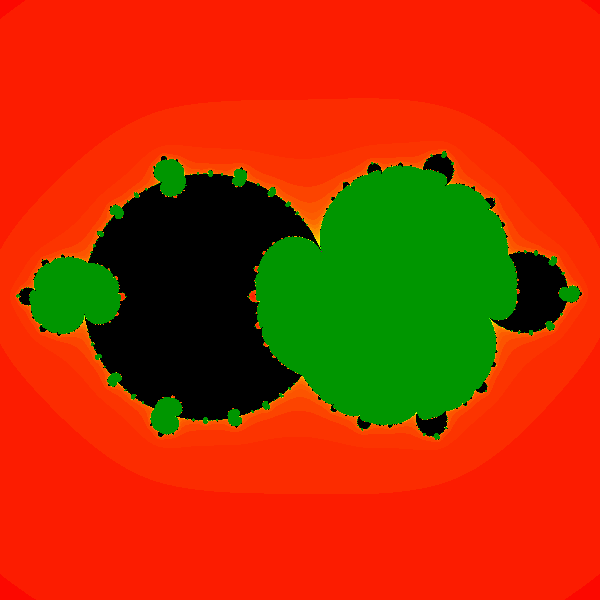}}
    \hspace{0.1in}  
 \subfigure[\scriptsize{Dynamical plane of $B_{2.55309+0.063042i}$.}  ]{
     \includegraphics[width=209pt]{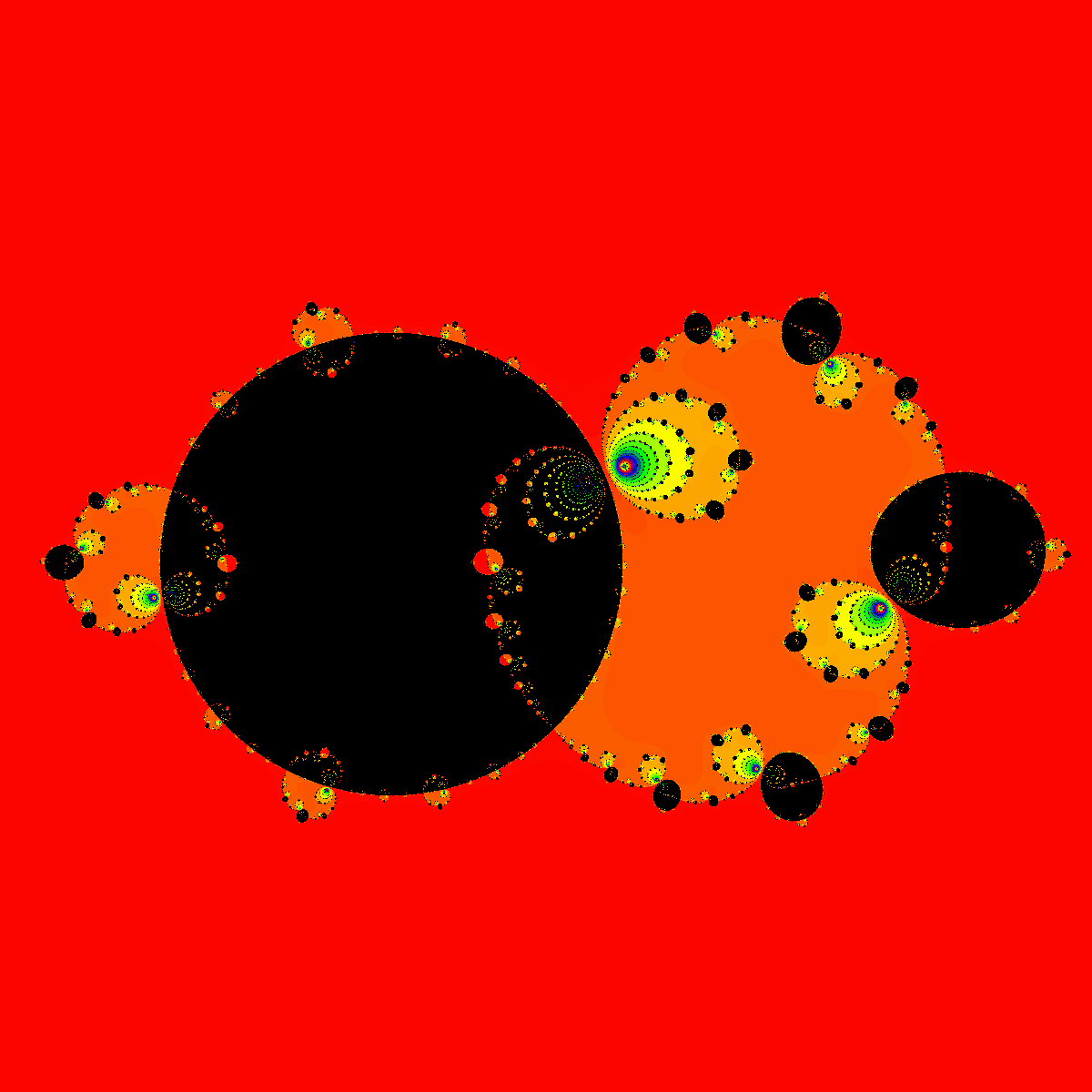}}
    \hspace{0.1in}  
   \subfigure[\scriptsize{Dynamical plane of $B_{2.64732+0.0421017i}$.}  ]{
    \includegraphics[width=209pt]{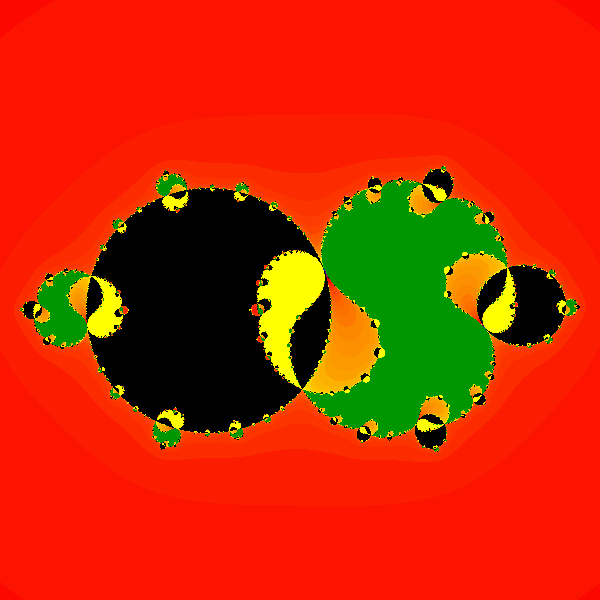}}
    \caption{\small{Figure (a) shows the dynamical plane of the tip $a_0=3$ of the fixed tongue $T_0$. Figure (b) shows the dynamical plane of $B_a$, where $a=2.65675+0.0389604i$ is in $\partial T_0\setminus a_0$. Figures (c) and (d) show the parameter plane of two Blaschke products with parameters near the boundary of the fixed tongue $T_0$. In Figure~(c) we have $a=2.55309+0.063042i$ and the parabolic fixed point has bifurcated into two repelling points while in Figure (d) we have $a=2.64732+0.0421017i$ and the parabolic fixed point has bifurcated into two attracting points. The colors are as follows: green if the point belongs to a basin of attraction which contains the critical point $c_+$, yellow if the point belongs to a basin of attraction which contains $c_-$ and not $c_+$, black if the orbit accumulates on $z=0$ and a scaling from blue to red if the orbit accumulates on $z=\infty$.  }} 
    \label{bifneartip}
\end{figure}

In this section we study the bifurcations which occur throughout the boundaries of the tongues. Given a tongue $T_{\tau}$, there is a persistent saddle-node bifurcation which takes place along $\partial T_{\tau}\setminus a_{\tau}$: two cycles collide in $\cercle$ and exit it (see Figure~\ref{bifneartip}).  The goal of the section is to prove Theorem~B, that is, to study the bifurcations in a neighbourhood of the tip of the tongues and to see that, if the parameter is close enough to a tip, the two cycles leaving the unit circle are attracting.

The proof of Theorem~B uses Proposition~\ref{finitestips}. Before stating the lemma, we introduce the concepts of algebraic geometry used on its proof. For a more detailed introduction to the topic we refer to \cite{Ha} and \cite{Sha}. First of all we recall the definition of algebraic variety in $\com^n$.
\begin{defi}
Let $\com\left[z_1,\cdots, z_n\right]$ be the polynomial ring of $n$ variables over $\com$. A subset $Y\subset \com^n$  is said to be an \textit{affine algebraic variety}  if there exists a finite set of  polynomials $S\subset \com\left[z_1,\cdots, z_n\right]$  such that $Y=\{ w\in\com^n\;| \; f(w)=0 \mbox{ for all } f\in S\} $.
\end{defi}

It is also useful to consider algebraic sets in the projective space $\mathbb{P}_{\com}^n$. In this setting we  work with homogeneous polynomials $f\in\com\left[z_0,\cdots, z_n\right]$ of degree $d$. Therefore,  if $a\in\com^{n+1}$ and $\lambda\in\com$, then $f(\lambda a)= \lambda^d f(a)$ and the sets of zeros on $\mathbb{P}_{\com}^n$ are well defined.  

\begin{defi}
 A subset $Y\subset \mathbb{P}_{\com}^n$  is said to be a \textit{projective algebraic variety}  if there exists a finite set of homogeneous  polynomials $S\subset \com\left[z_0,\cdots, z_n\right]$   such that $Y=\{ w\in\mathbb{P}_{\com}^n\;| \linebreak f(w)=0 \mbox{ for all } f\in S\} $.
\end{defi}

Affine varieties can be seen as subsets of projective varieties. Indeed, given polynomial $f\in \com\left[z_1,\cdots, z_n\right]$, it is not difficult to obtain a homogeneous polynomial $F\in \com\left[z_0,\cdots, z_n\right]$ such that if $w=(w_1,...,w_n)\in\com^n$ satisfies $f(w)=0$, then $w'=[1:w_1:\cdots:w_n]\in\mathbb{P}_{\com}^n$ satisfies $F(w')=0$.
 The Euclidean topology is not convenient when working with algebraic varieties. We shall work with the Zariski topology instead.

\begin{defi}
The \textit{Zariski topology} on $\com^n$ or $\mathbb{P}_{\com}^n$ is defined as the topology whose open sets are the complements of algebraic varieties.
\end{defi}

The Zariski topology is  well defined  since the finite union and intersection of algebraic varieties are algebraic varieties and the empty set and $\com^n$ (or $\mathbb{P}_{\com}^n$) are also algebraic varieties. Using this topology we can introduce the concepts of constructible sets and quasiprojective varieties.

\begin{defi}
A set $A$ is \textit{locally closed} if it is the intersection of an open and a closed subsets of $\com^n$ (or $\mathbb{P}_{\com}^n$). A set is said to be  \textit{constructible} if it is a finite disjoint union of locally closed sets.
\end{defi}

\begin{defi}
A \textit{quasiprojective variety} is a locally closed subset of an algebraic variety, i.e.\ the intersection of an algebraic variety with an open set of the Zariski topology.
\end{defi}

Notice that, by definition, quasiprojective varieties are constructible sets.
We finish this introduction to algebraic geometry presenting a simplified version of Chevalley's Theorem (see \cite{Gr}, c.f.\ \cite[Exercise 3.19]{Ha}), which is the main result we use on the proof of Proposition~\ref{finitestips}. 

\begin{thm}[Chevalley's Theorem]\label{Chevalley}
Any morphism of quasiprojective varieties sends constructible sets to constructible sets.
\end{thm}

We have not defined what a morphism of quasiprojective varieties is. However, we want to remark that we apply Chevalley's Theorem to the map given by the projection over one of the components of the variety, which is indeed a morphism of quasiprojective varieties. 

Now we can state and prove Proposition~\ref{finitestips}.

\begin{propo}\label{finitestips}
For fixed $n>0$, there is only a finite number of parameters $a\in\com$ for which the Blaschke product $B_a$ has a parabolic cycle of exact period $n$, multiplier 1 and multiplicity $3$.

\end{propo}
\proof
It will be convenient to work with the alternative parametrizations of the Blaschke products $B_a$ presented in Section~\ref{prelimparametrization}. Recall that, if $a=re^{2\pi i\alpha}$ with  $r>0$ and $\alpha\in\real$, $B_a$ is conjugate with $B_{r,3\alpha}$ (Equation (\ref{blasformula2})). The Blaschke products  $B_{r,3\alpha}$  are embedded within the family $G_{a,b}$ (Equation (\ref{Gab})), where $a,b\in\com$. We will prove that, for fixed $n>0$, there is only a finite number of parameters $(a,b)$, where $a,b\in\com$, for which  $G_{a,b}$ has a parabolic cycle of exact period $n$, multiplier 1 and multiplicity $3$.

We first show that the immediate basin of attraction of such a cycle contains both free critical points. Indeed, a parabolic cycle $\langle z_0\rangle$ of exact period $n$, multiplier~1 and multiplicity $3$ has two disjoint cycles of maximal attracting petals attached to it (see \cite[Theorem~10.7]{Mi1} and Figure~\ref{bifneartip} (a)). Each of these cycles of maximal petals has at least one critical point on the boundary of one of its components (see \cite[Theorem~10.15]{Mi1}). Therefore, the immediate basin of attraction of $\langle z_0\rangle$ contains both free critical points of $G_{a,b}$. Notice also that it follows from this assumption that the multiplicity cannot be greater than $3$ since the rational maps $G_{a,b}$ only have two free critical points. 

Parameters which satisfy the hypothesis are solutions of the system of rational equations

\begin{equation}\label{system0}
\left\{\begin{array}{l}
G_{a,b}^n(z)=z,\\
(\frac{\partial}{\partial z}G_{a,b}^n)(z)=1,\\
(\frac{\partial^2}{\partial z^2}G_{a,b}^n)(z)=0.\\
\end{array}\right. 
\end{equation}

Take the polynomials $p_1(z,a,b)$, $p_2(z,a,b)$, $p_3(z,a,b)$ and $q(z,a,b)$ so that the previous system reduces to

$$\left\{\begin{array}{l}
p_1(z,a,b)/q(z,a,b)=z,\\
p_2(z,a,b)/q(z,a,b)^2=1,\\
p_3(z,a,b)/q(z,a,b)^3=0.\\
\end{array}\right. $$

Notice that $p_2$ and $p_3$ are combinations of $p_1$, $q$ and their derivatives. We obtain the polynomial system of equations

\begin{equation}\label{system1}
\left\{\begin{array}{l}
p_1(z,a,b)-zq(z,a,b)=0,\\
p_2(z,a,b)-q(z,a,b)^2=0,\\
p_3(z,a,b)=0.\\
\end{array}\right. 
\end{equation}

The solutions of (\ref{system0}) also solve (\ref{system1}). However, we have added solutions. They come from points $(z,a,b)$ on which either the numerator and the denominator vanish simultaneously or are both equal to infinity. They can be equal to infinity if and only if $z=\infty$ or $a=\infty$ or $b=\infty$. Such points are not solutions of the original system. The point $z=\infty$ is a permanent superattracting fixed point (unless $b=0$ or $b=\infty$) and, therefore, does not satisfy the equations of a parabolic point. If $a=\infty$ then $G_{a,b}(z)$ degenerates to $bz^2$, which does not have parabolic cycles. If $b=\infty$  then $G_{a,b}$ is constant and therefore does not have any parabolic cycle. The points for which the numerator and the denominator vanish simultaneously come from the system

\begin{equation}\label{system2}
\left\{\begin{array}{l}
p_1(z,a,b)=0,\\
q(z,a,b)=0.
\end{array}\right. 
\end{equation}

Notice that, if $q(z,a,b)=0$ but $p_1(z,a,b)\neq 0$ then the first equation in (\ref{system1}) is not satisfied and, therefore the point $(z,a,b)$ is not a solution. Assume that $(z,a,b)$ solves (\ref{system2}). Then there is a $z_0$ such that the numerator and the denominator of

$$bG_{a,b}^{n-1}(z_0)^3\frac{G_{a,b}^{n-1}(z_0)-a}{1-aG_{a,b}^{n-1}(z_0)}$$

\noindent vanish simultaneously. This can only happen if $b=0$ or $a=\pm1$. If $b=0$ the map $G_{a,b}$ is constant and therefore  (\ref{system0}) has no solution.  If $a=\pm1$ the family $G_{a,b}(z)$ degenerates to the polynomials $\mp bz^3$ and the system (\ref{system0}) has no solution. 

We also assumed that the parabolic cycle has exact period $n$. Thus, the parameters which satisfy the hypothesis of the lemma are such that the equality 

\begin{equation}\label{system3}
G_{a,b}^{m}(z)=z
\end{equation}
\noindent is not satisfied for any $m<n$. If $G_{a,b}^{m}(z)=\tilde{p}_m(z,a,b)/\tilde{q}_m(z,a,b)$, the set of points which satisfy the previous equality, are solutions of the polynomial equation

\begin{equation}\label{system4}
\tilde{p}_m(z,a,b)-z\tilde{q}_m(z,a,b)=0.
\end{equation}

The set of solutions of (\ref{system1}) is an algebraic variety, say $Y$. Each point of $Y$ is either solution of (\ref{system0}) or corresponds to any of the degeneracy situations already described. The set of solutions of (\ref{system4}) consists of the solutions of (\ref{system3}) and exactly the same degeneracy solutions described for (\ref{system1}).  Let $Y'$ be the quasiprojective variety obtained by intersecting  $Y$ with the open set of the Zariski topology given by $b\neq 0$, $b\neq\infty$,  $a\neq\pm1$, $a\neq\infty$, $z\neq\infty$ and $\tilde{p}_m(z,a,b)-z\tilde{q}_m(z,a,b)\neq0$ for all $m<n$.  If $(z,a,b)$ belongs to $Y'$ then it is a solution of (\ref{system0}) and does not solve (\ref{system3}). Therefore, $\langle z\rangle$ is a parabolic cycle of $G_{a,b}$ of period exactly $n$, multiplier $1$ and multiplicity $3$ whose immediate basin of attraction contains both free critical orbits. Since by Lemma~\ref{boundednonescapingGab1} the non-escaping set is bounded in $a$, we conclude that the projection of $Y'$ over the variable $a$ is bounded.  It follows from Chevalley's Theorem (Theorem~\ref{Chevalley})  that the projection of a quasiprojective variety over a variable is a constructible set.  We conclude that the projection of $Y'$ over $a$ is finite since constructible sets in $\com$ are either dense in $\com$ or finite.  Summarizing we have that there are finitely many $a$ for which (\ref{system0}) has solution. 

Finally,  consider the previous equation systems with $a_0\neq\pm 1$ fixed. Let $Y'$ be the quasiprojective variety of points $(z,a_0,b)$ which solve (\ref{system0}), do not solve (\ref{system4}) for any $m<n$ and the degeneracy conditions are not satisfied. By Lemma~\ref{boundednonescapingGab2} we know that, for fixed $a_0\neq \pm1$, the non-escaping set is bounded on $b$. As before we conclude that $Y'$ projects onto a finite number of $b$, which finishes the proof.

\endproof

Notice that the condition of having exactly period $n$ on the previous proposition is necessary. Indeed, the family $B_a$ has curves of parabolic parameters  whose parabolic cycle $\langle z_0\rangle$ of period $n$ has multiplier $-1$ (see Theorem~C). The point $z_0$ is also a parabolic fixed point of $B_a^{2n}$ of multiplier $1$ and multiplicity $3$. Therefore, if we do not require exact period $n$ then we may obtain infinitely many solutions of (\ref{system0}).

\begin{figure}[hbt!]
\centering
\includegraphics[width= 11cm]{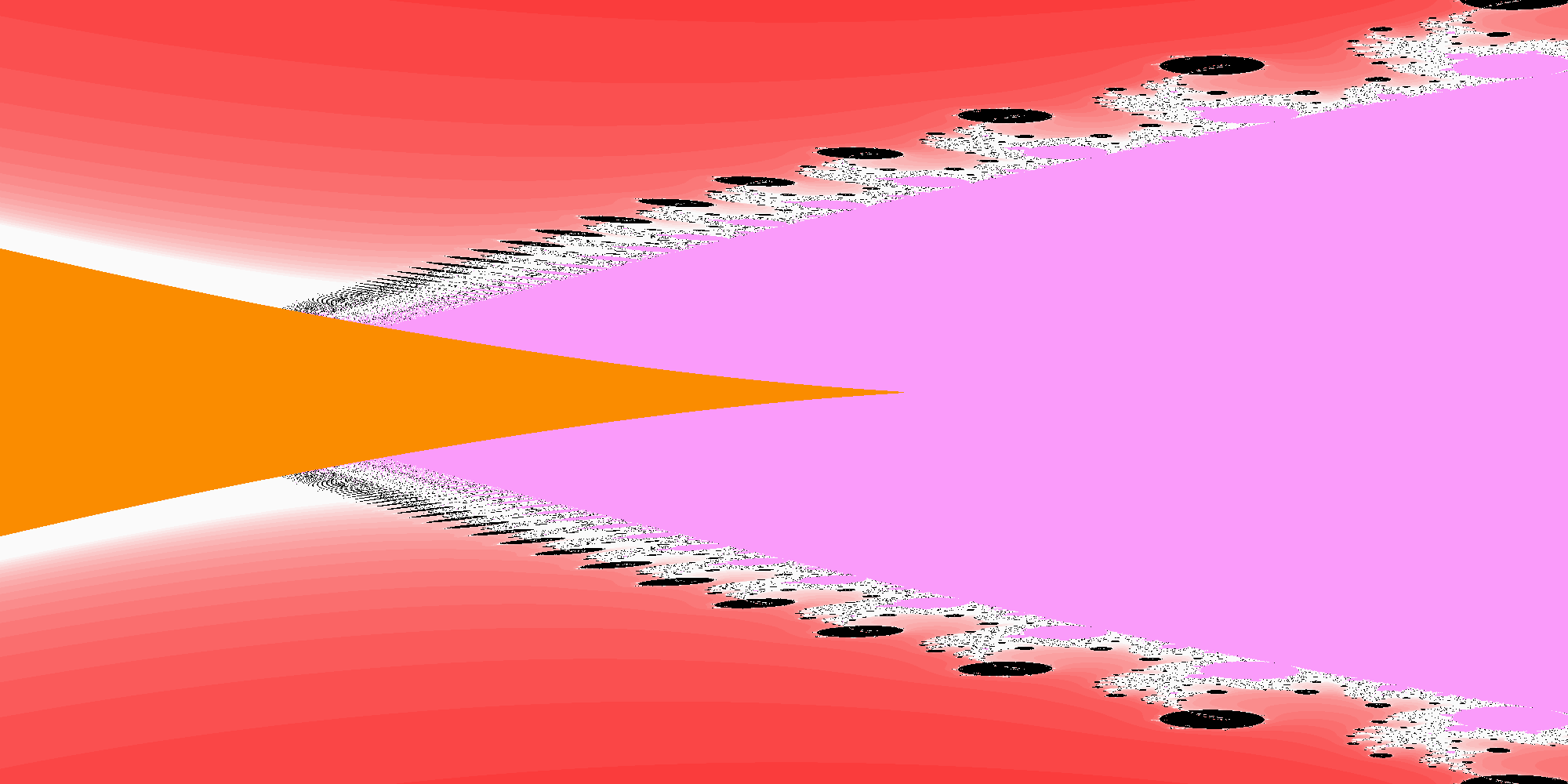}
  \put(-130,75){\color[rgb]{0,0,0}\makebox(0,0)[lb]{\smash{$a_0$}}}%
\caption{\small A zoom in a neighborhood of the tip $a_0$ of the tongue $T_o$. The colors are as in Figure~\ref{paramblash0}. We see in orange the fixed tongue $T_0$ and in pink a disjoint hyperbolic component which partially shares the boundary with $T_0$.}
\label{zoomtip}
\end{figure}

We now prove Theorem~B, which tells us that there is a neighborhood $U$ of the tip of any tongue such that if $a\in U$ then either $a$ belongs to the tongue, or to its boundary, or $B_a$ has two disjoint attracting cycles (see Figure~\ref{zoomtip} and Figure~\ref{bifneartip} (a), (b) and (d)).

\proof[Proof of Theorem~B]

The main ingredient for the proof is the holomorphic index. Given a fixed point $z_0$ of a holomorphic function $f$, the holomorphic fixed point index of $z_0$, denoted by $i(z_0)$, is defined to be the residue of $1/(z-f(z))$ around $z_0$. If  the fixed point has multiplier $\rho\neq 1$, then $i(z_0)=1/(1-\rho)$ (see \cite{Mi1}). Moreover, when $n$ different fixed points collide in a parabolic point $z_0$ of multiplier $1$, their indexes tend to infinity, even if the sum of their indexes tends to the finite index $i(z_0)$ of the parabolic point.

Let $\langle w_0\rangle$ be the parabolic cycle of $B_{a_{\tau}}$. Then, $w_0$ is a parabolic periodic point of multiplier 1, multiplicity $3$ and exact period $p$ of $B_{a_{\tau}}$.  Since, by Proposition~\ref{finitestips}, there is a finite number of such parameters, it follows that there is an open neighborhood of the parameter $a_{\tau}$ which contains no other parameter $a$ for which $B_a$ has a parabolic cycle of multiplier 1, multiplicity $3$ and the same period than $\langle w_0\rangle$.  Take a parameter $a$ close to $a_{\tau}$. The map $B_a^p$ has three fixed points, say $z_0$, $z_+$ and $z_-$, which tend to $w_0$ when $a$ tends to $a_{\tau}$. By symmetry and continuity of the semiconjugacy $H_a(x)$ with respect to $a$ and $x$, at least one of the fixed points lies in $\cercle$, say $z_0$, and satisfies $H_a(z_0)=\tau$. Also by symmetry, if more than one fixed point lies in $\cercle$, the three of them do.  In that later case, since $B_{a}|_{\cercle}^p$ is strictly increasing, one of them is either parabolic or attracting and satisfies $H_a(z)=\tau$ by continuity of $H_a$, so either belongs to the tongue $T_{\tau}$ or its boundary. Assume that only $z_0$ lies in $\cercle$ and is repelling (if it was attracting it would belong to $T_{\tau}$ again by continuity of $H_a$). Then $z_0$ has real multiplier $\eta>1$ (compare Lemma \ref{one attractingpar}). Due to the symmetry, the other two fixed points,  $z_{\pm}$, are symmetric. Moreover, their multipliers are complex conjugate say $\rho$ and $\overline{\rho}$. Indeed, we can conjugate $B_a$ via a M\"obius transformation $M$ to a rational map $\widetilde{B}_a$ that fixes the real line. The assertion follows then from the fact that $\widetilde{B}_a'(\bar{z})=\overline{\widetilde{B}'_a(z)}$ and that $M$ preserves the multiplier of the periodic cycles.

 Consider the sum $\mathcal{S}$ of the indexes of the three periodic  points.
$$\mathcal{S}=i(z_0)+i(z_+)+i(z_-)=\frac{1}{1-\eta}+\frac{1}{1-\rho}+\frac{1}{1-\overline{\rho}}.$$

The number $\mathcal{S}$ is a real quantity which tends to the index of the parabolic cycle of the tip of the tongue whenever $a$ tends to $a_{\tau}$. Moreover, $i(z_0)$ tends to minus infinity when $a$ tends to $a_{\tau}$. Hence, there is an open neighborhood $U$ of $a_0$ such that if $a\in U$ then $\tilde{\mathcal{S}}=\mathcal{S}-i(z_0)>1$. Write $\rho=1+\epsilon=1+\epsilon_r+i\epsilon_i$, $\epsilon_r, \epsilon_i\in \real$. Then, if $a\in U$, we have

$$\tilde{\mathcal{S}}=\frac{1}{1-1-\epsilon}+\frac{1}{1-1-\overline{\epsilon}}=-\frac{2\epsilon_r}{|\epsilon|^2}.$$

It follows from this equation that, if $a\in U$, then $\epsilon_r<0$. Finally, using that $|\epsilon|^2=-2\epsilon_r/ \tilde{\mathcal{S}}$, we have

$$ |\rho|^2=(1+\epsilon_r)^2+\epsilon_i^2=1+2\epsilon_r-\frac{2\epsilon_r}{\tilde{\mathcal{S}}} =1+2\epsilon_r(1-\frac{1}{\tilde{\mathcal{S}}}).$$

Since $\epsilon_r \lesssim 0$ and $\tilde{\mathcal{S}}>1$ we conclude that $|\rho|<1$, which finishes the proof.
\endproof

We finish the section showing some consequences of the construction presented in the proof of Theorem~B. The next corollary follows from the previous theorem and Theorem~\ref{simplygraph}, which states that the boundary of any tongue corresponds to the union of two arcs which intersect at the tip of the tongue. These two arcs can be parametrized univalently with respect to the modulus of the parameter.

\begin{co}
Given a tongue $T_{\tau}$, there exists a hyperbolic component of $B_a$ of disjoint type sharing part of its boundary with $T_{\tau}$ in a neighborhood of the tip $a_{\tau}$. 
\end{co}

The following proposition tells us that all parameters  $a$, with $|a|>2$, such that $B_a$ has a parabolic cycle which is topologically repelling in the unit circle are tips of tongues.

\begin{propo}\label{noparsol}
If $|a_0|>2$ and $B_{a_0}$ has a parabolic cycle of multiplicity $3$ and exact period $n$ in the unit circle then $a_0$ is the tip of a tongue of period $n$.
\end{propo}
\proof
Assume that $a_0$ is not a tip of a tongue of period $n$.  Then the same perturbation done in the proof of Theorem~B can be performed, obtaining a disjoint hyperbolic component of parameters which surrounds $a_0$. Indeed, since $a_0$ is not in the boundary of a tongue of period $n$ and there is a finite number of parameters with a parabolic cycle of multiplier 1, multiplicity $3$ and exact period $n$ by Proposition~\ref{finitestips}, the perturbation presented in the proof of Theorem~B gives us an open neighborhood $U$ of $a_0$ such that if $a\in U$, $a\neq a_0$, then the Blaschke product $B_a$ has two disjoint attracting cycles other than $z=0$ and $z=\infty$. Therefore, the set of parameters $U$ would be contained in a multiply connected disjoint hyperbolic component whose attracting cycles are not in the unit circle, which is impossible by Theorem~\ref{thmparametrize}.
\endproof

\begin{co}\label{parabolicinboundary}
If $|a|>2$ and $B_a$ has a parabolic cycle on the unit circle, then $a$ belongs to the boundary of a tongue.
\end{co}
\proof
If the parabolic cycle of $B_{a_0}$ has multiplicity $2$ then it is on the boundary of a period $n$ tongue by Lemma~\ref{clasification0}. If it has multiplicity $3$ it is on the tip of a fixed tongue by Proposition~\ref{noparsol}.

\endproof

\section{Extended Tongues}\label{exttongue}

The goal of this section is to give an idea of the dynamics that may take place for parameters within the open annulus $\mathbb{A}_{1,2}$ of inner radius 1 and outer radius 2. The section is structured as follows. We first notice that the tongues studied up to this point may be extended within this annulus. Then we describe more precisely how the fixed tongue extends. Finally we do some numerical computations to obtain an idea of other phenomena which may take place.

\begin{figure}[hbt!]
    \centering
    \subfigure{
    \includegraphics[width=209pt]{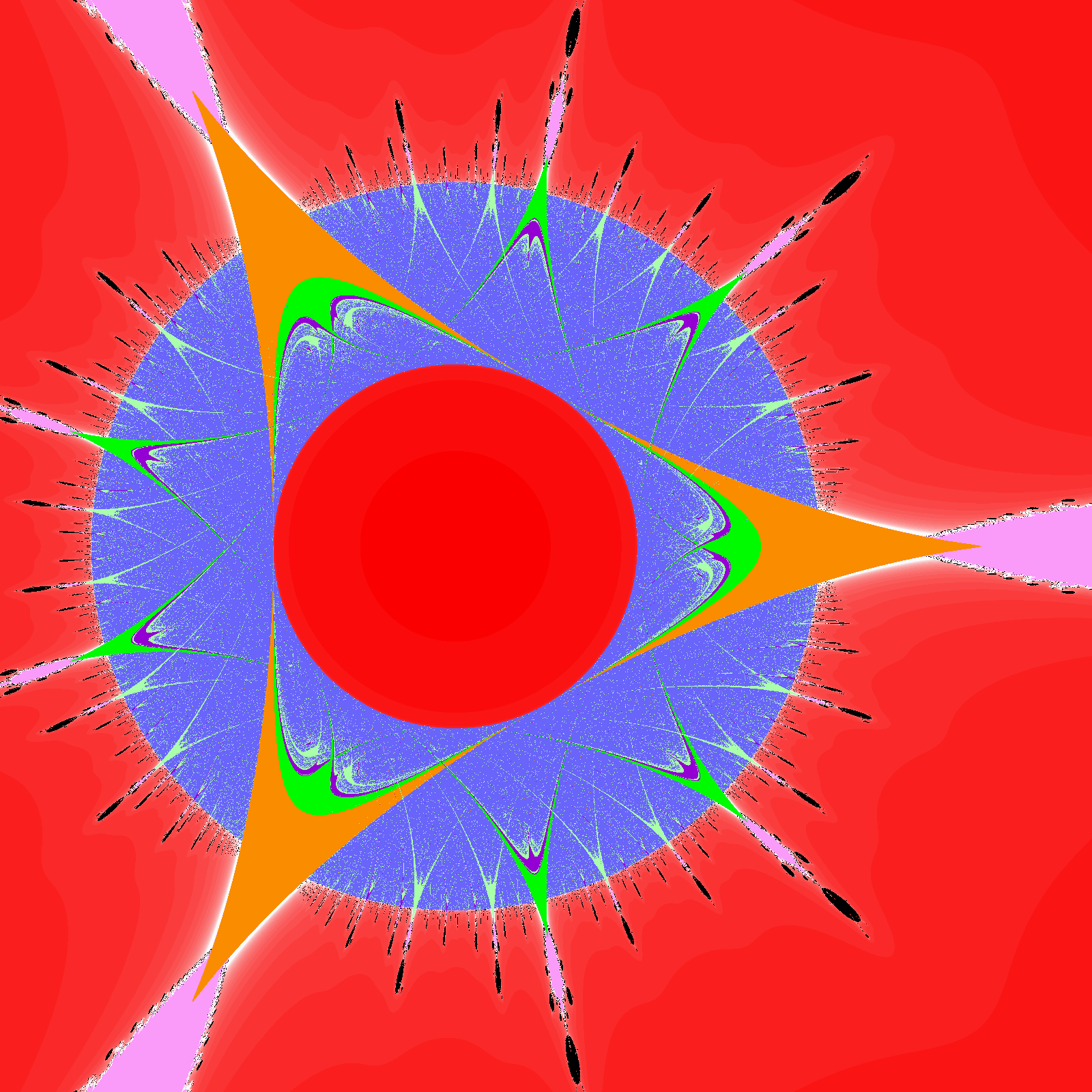}}
    \hspace{0.1in}
    \subfigure{
    \includegraphics[width=209pt]{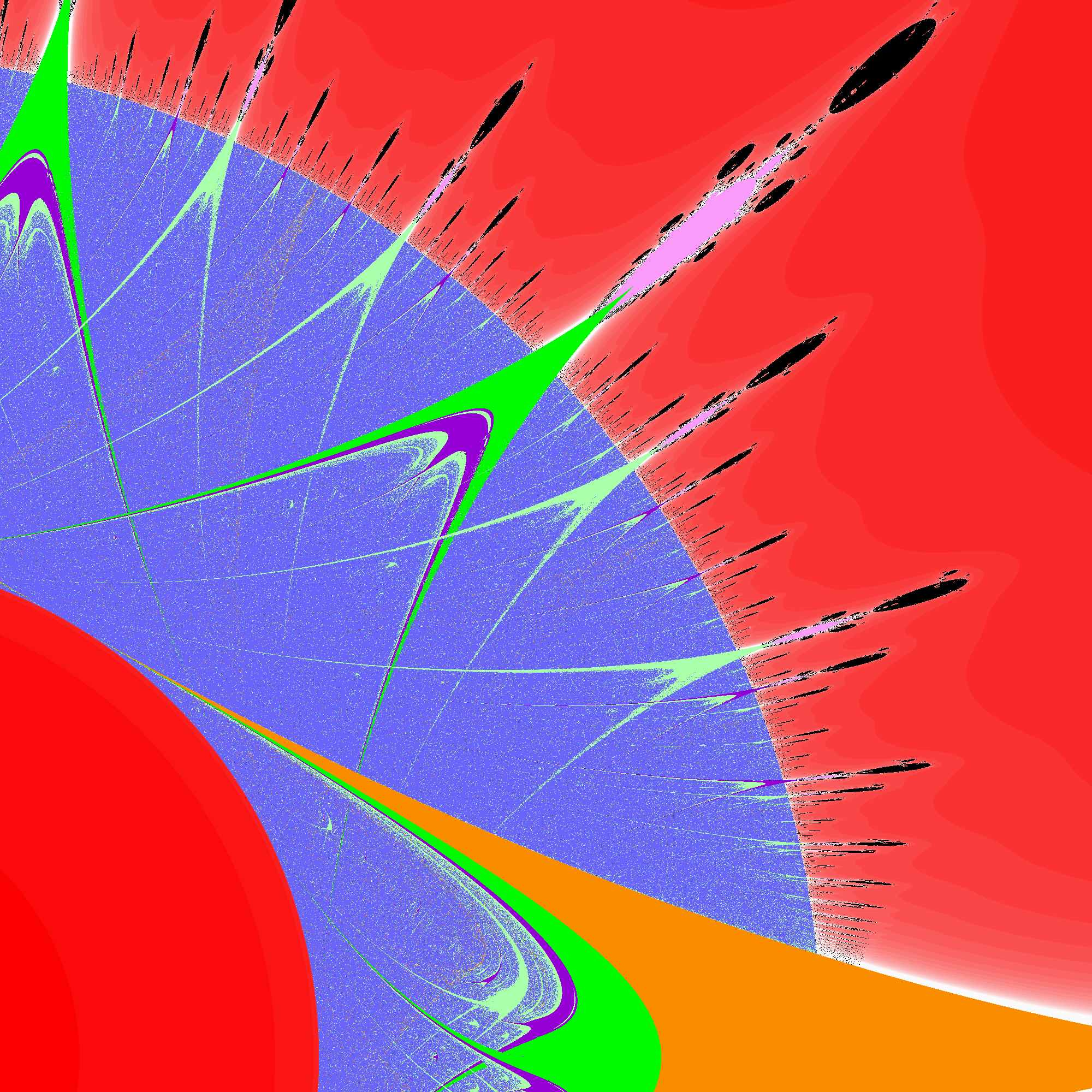}}
    
    \caption{\small Two zooms in Figure~\ref{paramblash0}. Recall that we plot in orange the parameters for which there is an attracting fixed point in $\cercle$. These parameters correspond to an extended fixed tongue. Strong green corresponds to parameters having a period 2 attracting cycle in  the unit circle, whereas violet corresponds to period 4 cycles. These parameters may belong to extended tongues of period 2 or 4,  or to other kinds of components.}
    \label{paramblash3col1}
\end{figure}

The definition of tongues only makes sense for parameters $a$ such that  $|a|\geq 2$. However, given a tongue $T_{\tau}$, its attracting cycle $\langle z_0\rangle$ can be analytically continued for parameters $1<|a|<2$ (see Figure~\ref{paramblash3col1}), parameters for which $B_{a}|_{\cercle}$  is not a degree 2 cover of the unit circle (see Section \ref{introblas0}) and, therefore, is not  semiconjugate to the doubling map.  We proceed to formally define the concept of extended tongue using the analytic continuation of the attracting cycle.

\begin{defin}
An \textit{extended tongue} $ET_{\tau}$ is defined to be the set of parameters for which the attracting cycle of $T_{\tau}$  can be continued analytically. More precisely, we say that a parameter $a$ belongs to the \textit{extended tongue} $ET_{\tau}$ (of period $p$)  if $1<|a|$, $B_{a}|_{\cercle}$ has an attracting periodic point of period $p$ and there exists a curve of parameters $\gamma(t)$ such that $\gamma(0)=a$, $\gamma(1)\in T_{\tau}$ and $B_{\gamma(t)}|_{\cercle}$ has an attracting fixed point of period $p$ for all $t\in(0,1)$ which depends continuously on $t$. 
\end{defin}

Since  the set of hyperbolic parameters is open in $\com$ and the roots of the tongues (parameters for which the cycle is superattracting) have modulus equal to 2 (see Theorem~A) we conclude that  $ET_{\tau}\cap\mathbb{A}_{1,2}$ is not empty for any periodic point $\tau$ of the doubling map. Notice that, with the previous definition, parameters $a\in T_{\tau}$ also belong to $ET_{\tau}$. The following lemma describes the parameters on the boundary of the extended tongues. Its proof is analogous to the one of Lemma~\ref{boundaryparabolic}.

\begin{lemma}
If $a$ belongs to the boundary of an extended tongue and $|a|\neq 1$, then $B_a$ has a parabolic periodic point of multiplier $\pm1$.
\end{lemma}

Notice that the intersection of two different extended tongues might be a non empty open set. Indeed, the critical orbits are not symmetric if $a\in\mathbb{A}_{1,2}$. Because of that, for $1<|a|<2$, $B_{a}|_{\cercle}$ might have two different attracting cycles, in which case $a$ might belong to two different tongues (see Figure~\ref{paramblash3col1}).

\subsection{The extended fixed tongue: proof of Theorem~C}

\begin{figure}[hbt!]
\centering
\includegraphics[width= 9cm]{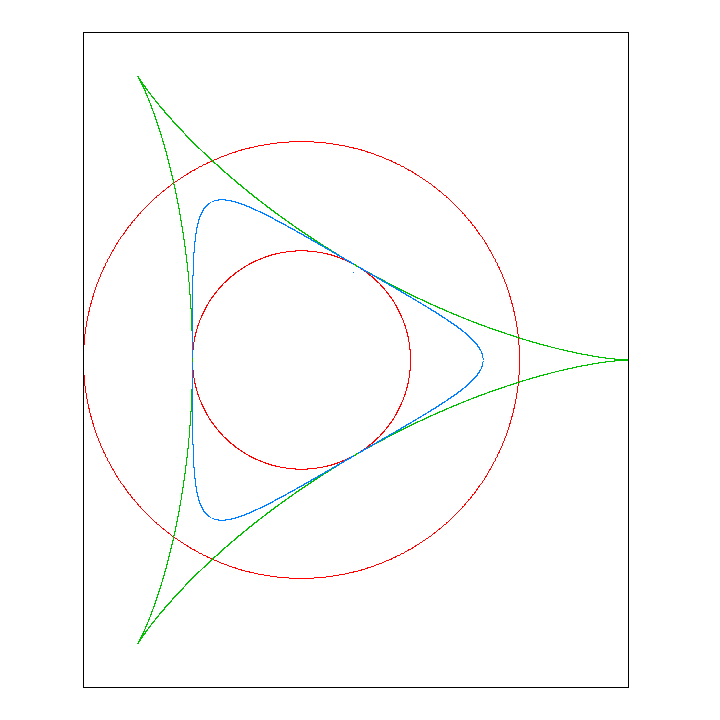}
\caption{\small Boundaries of the three symmetric extended fixed tongues. The green curves correspond to parameters for which a fixed point has multiplier 1. The blue curves correspond to parameters for which a fixed point has multiplier -1. 
 We also plot in red the two circles of parameters $|a|=1$ and $|a|=2$.}  
\label{corba}
\end{figure}

The goal of this subsection is to proof Theorem~C, which describes the shape of the connected components of the extended fixed tongue $ET_0$ (see Figure~\ref{corba}).  

We need an auxiliary proposition to prove Theorem~C. The fixed tongue $T_0$ has three connected components, only one modulo symmetry (see Theorem~A). Therefore, when studying the extended fixed tongue we restrict to the connected component which intersects the real line.   It is convenient to consider the parameter plane given by $(r,\alpha)$, where $a=r e^{2\pi i \alpha}$,  $1<r<2$ and $\alpha\in [0, 1/3)$. We then use the alternative parametrization $g_{r,\alpha}=B_{r,3\alpha}|_{\cercle}$ of $B_{a}|_{\cercle}$ (see Section~\ref{prelimparametrization}). We denote by $h_{r,\alpha}$ the lift of $g_{r,\alpha}$ (see Equation~\ref{hra}).  We want to remark that for $r=1$ and $x=0$ the function is not well defined (see Section~\ref{introblas0}). Indeed, for $r=1$, the two critical points and the preimages of $0$ and $\infty$ collapse at the point $x=0$ and the function becomes a degree $3$ polynomial.  The following proposition gives us the main properties of $ET_0$.

\begin{propo}\label{extendedfixed}
Let $ET_0$ denote the extended fixed tongue which intersects the real line. Then, $ET_0$ satisfies the following properties:
\begin{enumerate}[(a)]
\item $ET_0$ is symmetric with respect to the real line.
\item For fixed $r_0$, $1<r_0<2$, $ET_0\cap \{\alpha\geq 0\}\cap \{r=r_0\}$ is a connected set on which the multiplier is strictly increasing with respect to $\alpha$ and takes values in $(b, 1)$, where $-1\leq b<0$. Moreover, $b=-1$ if and only if $r_0\leq 5/3$.

\item If $(r,\alpha)\in ET_0$, then  $-1/6<\alpha<1/6$.
\end{enumerate}

\end{propo}

\proof

For statement (a) we notice  that the extended tongue is symmetric with respect to the real line due to the symmetry of the parameter plane $a\rightarrow \bar{a}$ (see Lemma \ref{symm}), so we restrict the study to $\alpha\geq 0$. We now use the lift $h_{r,\alpha}$ of $g_{r,\alpha}$ (see Equation~(\ref{hra})) to prove (b) and (c). Recall that it is given by 

$$h_{r,\alpha}(x)=3x+3\alpha+\frac{1}{2\pi i}\log\left(\frac{e^{2\pi i x}-r}{1-re^{2\pi i x}}\right).$$

 The extended fixed tongue $ET_0$ consists of the parameters for which $h_{r,\alpha}$ has an attracting fixed point which can be continued to the superattracting fixed point $x_{r,0}=0$ of $h_{2,0}$, where the parameter $(2,0)$ corresponds to the root $2=2e^{2\pi i 0}$ of the fixed tongue $T_0$. 

The point $x_{r,0}=0$ is a fixed point of $h_{r,0}$ for all $r\in (1,2]$. We use the monotonicity of $h_{r,\alpha}$ with respect to $\alpha$ to continue this fixed point for $\alpha>0$. We denote by $x_{r,\alpha}$ the fixed point of $h_{r,\alpha}$ obtained by continuation of $x_{r,0}=0$. This point $x_{r,\alpha}$ is well defined as long as we do not reach a parameter $\alpha$ for which $x_{r,\alpha}$ has multiplier $1$. Such a parameter would belong to the boundary of the extended tongue, so this is not an obstruction. Notice also that $x_{r,\alpha}$ is strictly decreasing with respect to $\alpha$ since $h_{r,\alpha}$ is strictly increasing with respect to $\alpha$.

To determine if a parameter $(r,\alpha)$ belongs to the extended fixed tongue, we have to study the multipliers of the fixed points $x_{r,\alpha}$. Recall from Equation~(\ref{derhra}) that the derivative of $h_{r,\alpha}$ is given by 

$$h'_r(x):=\frac{\partial}{\partial x} h_{r,\alpha}(x)=3+\frac{1-r^2}{1+r^2-2r \cos(2\pi x) }.$$

Notice that it does not depend on $\alpha$. The multiplier of the fixed point $x_{r,0}$ is given by $h'_r(0)=3+(1+r)/(1-r)$,  which decreases from $0$ to $-\infty$ as $r$ decreases from $2$ to $1$. Moreover, $h'_r(5/3)=-1$. The proof of (b) will be finished finding, for every fixed $r$ with $1<r\leq2$, a parameter $\alpha_1(r)$ such that $x_{r,\alpha_1(r)}$ has multiplier $1$ and such that the multipliers of $x_{r,\alpha}$ increase from $h'_r(0)$ to 1 as $\alpha$ increases from 0 to $\alpha_1(r)$.

Fixed $r$ with $1<r\leq2$, the monotonicity of the multiplier of $x_{r,\alpha}$ with respect to $\alpha$ follows from the fact that $h'_r$ is strictly increasing as $x$ decreases from $0$ to $-1/2$. The condition $0>x>-1/2$ is not an obstruction since $h'_r(-1/2)=3+(1-r^2)/(1+r^2+2r)>1$ for all $1<r\leq 2$ and, therefore, no parameter $(r, \alpha)$ with $x_{r,\alpha}=-1/2$ can belong to the extended fixed tongue or its boundary.

We now prove that for all $r\in(1,2]$, there exists an $\alpha_1(r)=\alpha_1$ which depends continuously on $r$ such that $h'_r(x_{r,\alpha_1})=1$ and $0<\alpha_1<1/6$. First of all we notice that, given a parameter $(r, \alpha_1(r))$ with a fixed point $x_{r, \alpha_1(r)}=x_1$ of multiplier 1, the map $h_{r,\alpha_1(r)}$ can be written as

$$h_{r,\alpha_1(r)}(x)=x+\eta(x-x_{1})^2+\mathcal{O}((x-x_{1})^3)$$

\noindent in a neighborhood of $x_1$, where $\eta\in\real$. We conclude that $\alpha_1(r)$ is a local graph with respect to $r$ unless $\eta=0$. However, $\eta$ is zero if $\partial^2h_{r}/\partial x^2 (x_1)=0$, which can only happen if $x_1=0$ or $x_1=-1/2$. Since  $h'_r(-1/2)=3+(1-r^2)/(1+r^2+2r)>1$,  $x=-1/2$ cannot be a parabolic fixed point. Moreover,  the point $x=0$ cannot be a parabolic fixed point with multiplier 1  since $h'_r(0)\leq 0$ if $1<r\leq2$. 

Next we prove that $\alpha_1(r)<1/6$ by contradiction. Assume that there is an $r$ for which this is not the case. Then, by continuity, there would be an $\tilde{r}$ so that $h'_r(x_{\tilde{r},1/6})=1$. Because of the symmetries in the parameter plane, this parameter would give us the intersection between the boundaries of the extension of two different connected components of the fixed tongue. Therefore, at this parameter we would have two different fixed points of multiplier 1. However, each of these parabolic points is a fixed point of multiplicity two and, therefore, this situation would require at least 4 fixed points. This would contradict the fact that the Blaschke products $B_a$ can have at most 3 fixed points other than $z=0$ and $z=\infty$. 

Summarizing we have proven that, for all $r\in(1,2)$, the fixed point $x_{r,0}=0$ of $h_{r,0}$ can be monotonously continued to a fixed point $x_{r,\alpha}$ of $h_{r,\alpha}$ as long as $\alpha<\alpha_1(r)<1/6$, where $\alpha_1(r)$ is a continuous function with respect to $r$. Moreover, for $0\leq\alpha<\alpha_1(r)$ the multiplier of $x_{r,\alpha}$ is strictly increasing and takes values in $[h'_r(0),1)$, where $h'_r(0)\leq-1$ if and only if $r\leq 5/3$. This finishes the proof of the proposition.

\endproof

Using the previous proposition we can prove Theorem~C.

\proof[Proof of Theorem~C]

The fact that two extensions of connected components of the extended fixed tongue cannot intersect follows from statement (c) of Proposition~\ref{extendedfixed} using the symmetries of the parameter plane and that all connected components of the fixed tongue are symmetric with respect to the rotations  given by the third roots of the unity (see Theorem~A).

The boundary of the connected components of $ET_0$ is the union of a exterior boundary with parameters of multiplier $1$ and a interior boundary with parameters of multiplier $-1$ by statement (b) of Proposition~\ref{extendedfixed}.

We finally prove that there is a period doubling bifurcation taking place throughout the interior curve. Let $(r,\alpha_{-1})$, $r<5/3$, be a parameter of the interior curve. Then $h_{r,\alpha_{-1}}$ has a parabolic fixed point of multiplier $-1$, say $x_-$. Hence, $x_-$ is a parabolic fixed point of multiplier $1$ of $h_{r,\alpha_{-1}}^2$. Using that $h_{r,\alpha_{-1}}(x_-)=x_-$ and $h'_{r,\alpha_{-1}}(x_-)=-1$ it is not difficult to prove that $\partial^2h_{r,\alpha_{-1}}^2/\partial x^2 (x_-)=0$ and, therefore, $x_-$ is a parabolic fixed point of multiplicity $3$ of $h_{r,\alpha_{-1}}^2$. Hence, $x_-$ has two attracting petals which intersect the unit circle since all critical points lie in $\cercle$. Consequently, $x_-$ is topologically attracting on the unit circle. Using that  $h_{r,\alpha_{-1}}^2$ is monotonously decreasing with respect to $\alpha$ in an open neighborhood of $x_-$ and performing the same perturbations as in Lemma~\ref{clasification0} we conclude that, if $\alpha\in(\alpha_{-1}-\epsilon,\alpha_{-1})\cup(\alpha_{-1},\alpha_{-1}+\epsilon)$ with $\epsilon>0 $ small enough, then $h_{r,\alpha}^2$ has a topologically attracting fixed point. It follows from the monotonicity of the multipliers of the fixed points of $h_{r,\alpha}$ with respect to $\alpha$ shown in Proposition~\ref{extendedfixed} that either the parameters $(r,\alpha)$ with $\alpha<\alpha_{-1}$ or the $(r,\alpha)$ with $\alpha>\alpha_{-1}$ are such that $h_{r,\alpha}$ has a period two attracting cycle.
\endproof

\subsection{Numerical conjectures}

We finish this section giving some ideas about the dynamics that take place on $\mathbb{A}_{1,2}$ other than the ones given by the extended fixed tongue. Numerical studies suggest that the properties described for the extended fixed tongue $ET_0$ are common to all extended tongues. Indeed, we conjecture that all extended tongues have a similar structure than the one  presented in Theorem~\ref{extendedfixed}.

\begin{conj}
Given an extended tongue $ET_{\tau}$ of period $p>1$, its connected components are disjoint. The boundary of every connected component of the extended tongue $ET_{\tau}$ consists of two disjoint connected components. The exterior component consists of parameters for which there is a parabolic cycle of period $p$ and multiplier $1$. The interior component consists of parameters for which there is a parabolic cycle of period $p$ and multiplier $-1$. Moreover, there is a period doubling bifurcation taking place throughout the curve of interior boundary parameters.
\end{conj}

We also  focus on how extended tongues, and more generally hyperbolic components, accumulate on the unit circle. In Figure~\ref{corba} we see that the boundaries of the extended fixed tongue accumulate tangentially to the unit circle onto isolate points. Moreover they do it in pairs. By this we mean that given an accumulation point $l$, $|l|=1$, the boundaries of two different connected components of the extended fixed tongue land on it. This seems to happen for hyperbolic regions of arbitrary period. In Figure~\ref{paramblash3col1} (right)) we can observe how the boundary of an extended tongue of period 2 (drawn in green) accumulates onto a parameter $l$, $|l|=1$, together with another hyperbolic region of period 2 not coming from the extension of a tongue.

\begin{conj}
Parabolic curves accumulate on the unit circle on isolate points and tangentially to it. Moreover, if a parabolic curve accumulates on a parameter $l$, $|l|=1$, then there is another parabolic curve of the same period landing on $l$ from the opposite side.
\end{conj}

\begin{figure}[hbt!]
\centering
\includegraphics[width= 6.5cm]{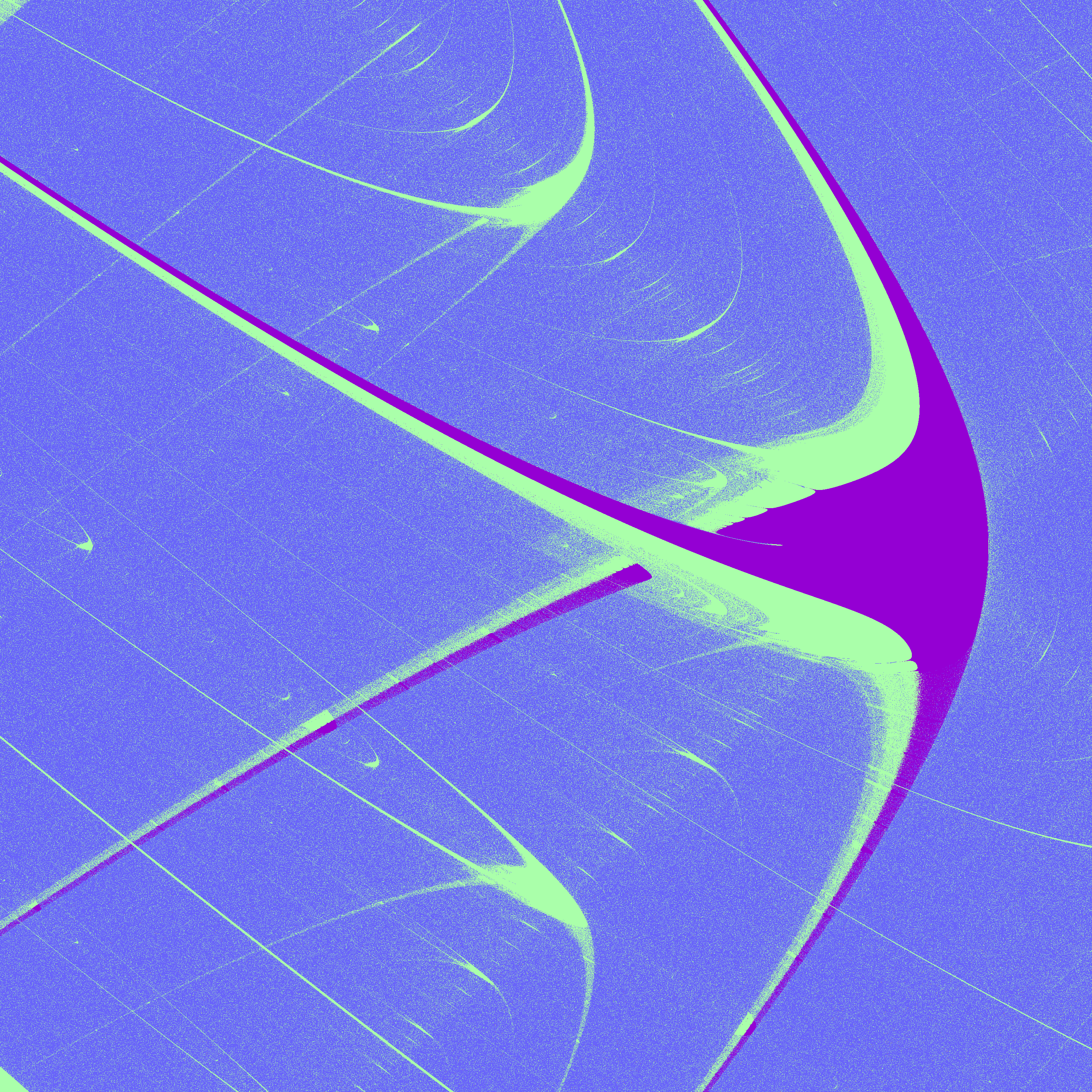}
\caption{\small Zoom in Figure \ref{paramblash0}. A period 4 Cross-road can be observed in violet. Parameters $a$ are taken so that $1.2<\re(a)<1.24$ and $-0.02<\im(a)<0.02$.}
\label{crossroad}
\end{figure}

 We also want to point out that the observed bifurcation structures are similar to the ones described for more general maps of $\real^2$ (see \cite{BST, CMBST, GSV}). Indeed, the observed structure around extended tongues is very similar to the one of spring-areas associated to homoclinic tangencies. In that later case, it is proven that there are cascades of period doubling bifurcations which are also observed for the family $B_a$. In Figure \ref{paramblash3col1} we can see these period doublings to period 2 and 4.  Other structures of bifurcations described in these papers such as Cross-roads also seem to appear for the Blaschke family $B_a$ when $1<|a|<2$ (see Figure \ref{crossroad}).

\bibliography{bibliografia}
\bibliographystyle{amsalpha}

\end{document}